\documentclass{svjour3}
\usepackage[utf8]{inputenc}
\usepackage{amsmath,amssymb}
\usepackage{hyperref}
\usepackage{graphicx}
\newcommand{\RR}{\mathbb R}
\newcommand{\SSS}{\mathbb S}
\newcommand{\CC}{\mathbb C}
\newcommand{\PP}{\mathbb P}
\newcommand{\QQ}{\mathbb Q}
\newcommand{\ZZ}{\mathbb Z}
\newcommand{\im}{\mathrm{im}\,}

\title{Bad Projections of the PSD Cone}
\author{Yuhan Jiang \and Bernd Sturmfels}
\institute{Yuhan Jiang 
\at Harvard University, MA 02138, USA, 
\\ \email{yjiang@math.harvard.edu} 
\and Bernd Sturmfels 
\at Max Planck Institute for Mathematics in the Sciences, Leipzig 04103, Germany, 
\\ \email{bernd@mis.mpg.de} 
}
\date{December 2020}

\begin{document}

\maketitle
\begin{abstract}
    The image of the cone of positive semidefinite matrices under a linear map is a convex cone.
    Pataki characterized the set of linear maps for which that image is not closed.
    The Zariski closure of this set is a hypersurface in the Grassmannian.
    Its components are the coisotropic hypersurfaces of symmetric determinantal varieties.
    We develop the convex algebraic geometry of such bad projections, with focus on explicit computations.
\end{abstract}

\section{Introduction}
\label{sec1}
Real symmetric $n \times n$ matrices are identified with quadratic forms on $\RR^n$, and they form a vector space $\SSS^n$ of dimension $\binom{n+1}{2}$. 
We write $\SSS^n_+$ for the subset of quadratic forms that are nonnegative on $\RR^n$. 
This is a full-dimensional closed semialgebraic convex cone  
in $\SSS^n$, known as the \emph{PSD cone}. 
Its elements are identified with positive semidefinite matrices.
The PSD cone is self-dual with respect to the trace inner product
$A \circ B  := \mathrm{trace}(A B)$ for $A, B \in \SSS^n$.

Given any linear subspace $\mathcal{L} $ of $\SSS^n$, we consider the linear projection $\pi_\mathcal{L}: \SSS^n \rightarrow \mathcal{L}^\vee$ that is dual to the inclusion $\mathcal{L} \subset \SSS^n$. 
Here, $\mathcal{L}^\vee = {\rm Hom}(\mathcal{L},\RR)$ denotes the vector space dual of $\mathcal{L}$.
We are interested in the image $ \pi_\mathcal{L}(\SSS^n_+)$
of the PSD cone under this map.
These objects can be written in coordinates as follows.
If  $\{A_1,A_2,\ldots,A_k\} \subset \SSS^n$
is a basis of $\mathcal{L}$ then our map is
\begin{equation}
\label{eq:pimap} \pi_\mathcal{L} \,:\, \SSS^n \,\to\, \RR^k \,, \,\, X \,\mapsto \,\bigl( A_1 \circ X,\,
 A_2 \circ X , \ldots,  \,A_k \circ X  \bigr). \end{equation}
While the PSD cone $\SSS^n_+$ is closed in $\SSS^n$, its image under the map $ \pi_\mathcal{L}$ may not be.
\begin{example} \label{ex:badintroex1}
Fix $n=k=2$. Let $\mathcal{L} $ be the linear space spanned by the two quadratic forms $q_1 = x_1^2$ and $q_2 = 2x_1 x_2$.
The corresponding symmetric $2 \times 2$ matrices  are
$$\begin{small} A_1 \,= \,\begin{pmatrix} 1 & 0 \\ 0 & 0 \end{pmatrix} \quad {\rm and} \quad
A_2 \, = \,  \begin{pmatrix} 0 & 1 \\ 1 & 0 \end{pmatrix} . \end{small} $$
The linear map $\pi_\mathcal{L}$ projects the $3$-dimensional space  $\SSS^2$ into the plane $\RR^2$ via
$$X = \begin{small} \begin{pmatrix} x_{11} & x_{12} \\ x_{12}  & x_{22} \end{pmatrix} \end{small}\, \,\mapsto \, \,
( A_1 \circ X,\, A_2 \circ X ) \, = \,
( \,x_{11}, \,2 x_{12} \,).$$
The image of the PSD cone $\SSS^2_+$ under this map into $\RR^2$ is not closed. We find that
$$\ \pi_\mathcal{L} ( \SSS^2_+) \,\,=\, \, \bigl\{ \,(z_1,z_2) \in \RR^2 \,: z_1 > 0
\,\,\,{\rm or}\,\,\, z_1 = z_2 = 0\, \bigr\}.$$
\end{example}
The failure of the image to be closed reflects the fact that strong duality can fail in semidefinite programming (SDP).
A thorough study of this phenomenon was undertaken by Pataki in \cite{LP15}, \cite{Pataki13}, \cite{Pataki17}, \cite{Pataki19}.
Basics on SDP and its algebraic aspects can be found in \cite[Chapter 1]{BPT} and \cite[Chapter 12]{MSbook}.
Our subspace $\mathcal{L}$ plays the role of an instance of SDP, as in \cite[Corollary 12.12]{MSbook}.
Using the adjective proposed in \cite{Pataki17}, the subspace $\mathcal{L}$ is called \emph{bad} if 
the cone $\pi_\mathcal{L} ( \SSS^n_+) $ is not closed.
But, just like in slang usage, ``bad'' can also mean ``good''.
Pataki derives a  characterization, and he concludes that \emph{bad semidefinite programs all look the same}~\cite{Pataki17}.

The aim of this paper is to examine this phenomenon through the lens of algebraic geometry. 
The assertion of \cite{Pataki17} that all bad instances ``look the same'' refers to the natural action of the group ${\rm GL}(n) \times {\rm GL}(k)$ on the domain and range of our map $\pi_\mathcal{L}$ in (\ref{eq:pimap}).
Pataki describes normal forms of bad instances $\mathcal{L}$ with respect to that action, to be reviewed in Section \ref{sec2}.
We are here interested in the geometry of the locus of all bad instances, that is, the orbits of Pataki's normal forms under ${\rm GL}(n) \times {\rm GL}(k)$.
We pass to the Zariski closure, and study the corresponding complex projective variety.
The following example is meant to illustrate how our perspective builds on and differs from that developed in \cite{LP15}, 
\cite{Pataki13}, \cite{Pataki17} and \cite{Pataki19}.

\begin{example} \label{ex:badintroex2}
Fix $n=k=2$. Let $\mathcal{L}$ be a general subspace in $\SSS^2$,
with a given basis 
$$ q_1 \,=\, a_{11} x_1^2 + 2 a_{12} x_1 x_2 + a_{13} x_2^2 
\quad {\rm and} \quad
q_2 \,=\, b_{11} x_2^2 + 2 b_{12} x_1 x_2 + b_{13} x_2^2.
$$
The nature of  $\pi_\mathcal{L}(\SSS_+^2)$ is determined by the \emph{resultant} of these two binary quadrics:
\begin{small}
$$ R \,= \,
a_{11}^2 b_{22}^2-4 a_{11} a_{12} b_{12} b_{22}-2 a_{11} a_{22} b_{11} b_{22}
+4 a_{11} a_{22} b_{12}^2 + 4 a_{12}^2 b_{11} b_{22}
-4 a_{12} a_{22} b_{11} b_{12}+a_{22}^2 b_{11}^2.
$$
\end{small}
Our theory in Section \ref{sec3} implies that \emph{$\pi_\mathcal{L}(\SSS_+^2)$ is not closed in $\RR^2$ if and only if $R=0$.}
Furthermore, if $R > 0$ then $\pi(\SSS_+^2) $ is a closed pointed cone. 
Yet, if $R < 0$ then $\,\pi(\SSS_+^2)  = \RR^2$.
If $ R > 0$ then $\mathcal{L}$ is spanned by two squares,
$\gamma_1 q_1 + \gamma_2 q_2 = (u_1 x_1 + u_2 x_2)^2$ and
$\delta_1 q_1 + \delta_2 q_2 = (v_1 x_1 + v_1 x_2)^2$, and we have
$\pi(\SSS^2_+) \,= \, \bigl\{ (z_1,z_2) \in \RR^2 \,:\,
\gamma_1 z_1 + \gamma_2 z_2 \geq 0  \,\,\,{\rm and} \,\,\,
\delta_1 z_1 + \delta_2 z_2 \geq 0  \,\bigr\}$.
If $R=0$ then the two squares are linearly dependent and the
image cone is not closed:
\begin{equation}
\label{eq:nonclosed}
 \pi(\SSS^2_+) \quad = \quad \bigl\{ \,(z_1,z_2) \in \RR^2 \,:\,
 \gamma_1 z_1 + \gamma_2 z_2 > 0 \,\bigr\} \,\,\cup \,\, \bigl\{ (0,0) \bigr\}.
\end{equation}
In conclusion, all bad instances do look the same as Example~\ref{ex:badintroex1}.
But, from an algebraic perspective, their parameter space $\{R=0\}$ is a variety of great interest.
\end{example}

This article is organized as follows. 
Section \ref{sec2}  characterizes all linear subspaces~$\mathcal{L}$ in $\SSS^n$ that are bad in the sense that the image cone $\pi_\mathcal{L} (\SSS^n_+)$ is not closed.
This result is due to Pataki \cite{Pataki17}, \cite{Pataki19}. 
In Theorem \ref{conditions} we present his linear algebra formulation in terms of block matrices. 
We recast this in the setting of real algebraic geometry. This is motivated by Proposition~\ref{prop:motivation} which states that the bad subspaces $\mathcal{L}$ form a semialgebraic subset of the Grassmannian ${\rm Gr}(k,\SSS^n)$.
In Example \ref{ex:disagreements}, we offer a contrast to the analogous closure question for the images of the quadratic maps $\RR^n \rightarrow \RR^k$ and $\CC^n \rightarrow \CC^k$ that also arise from our subspace $\mathcal{L} \simeq \RR^k$ of $\SSS^n$.

In Section \ref{sec3}, we turn to projective geometry and study the Zariski closure ${\rm Bad}_{k,n}$ of the set of bad subspaces $\mathcal{L}$ inside the complex Grassmannian ${\rm Gr}(k,\SSS^n)$. 
These varieties have codimension one and they are generally reducible.
Their irreducible components are the coisotropic hypersurfaces \cite{kohn} of rank strata of symmetric matrices. 
This is the content of Theorem \ref{thm:BadCoiso}, which identifies bad projections with objects familiar from elimination theory, such as resultants, Chow forms and Hurwitz forms \cite{GKZ}, \cite{Stu}.
Hyperdeterminants \cite{GKZ} remain in the background.
Examples \ref{ex:nn33} and \ref{ex:vier} offer detailed analyses of the cases $n=3$ and~$n=4$.

Section \ref{sec4} explains the badness of a subspace $\mathcal{L}$ in terms of the normal cycle of the cone $\SSS^n_+$ and its Zariski closure in $\PP(\SSS^n) \times \PP(\SSS^n)$. 
The latter is the projective normal cycle, whose irreducible components are the conormal varieties of the rank strata \cite[Example 5.15]{BPT}. 
Following \cite{NRS}, this encodes complementarity in SDP.
Theorem \ref{thm:LtimesLperp} reveals that $\mathcal{L}$ is bad when the normal cycle meets $\mathcal{L} \times \mathcal{L}^\perp$.
This furnishes effective algebraic tools to identify bad projections, illustrated by computations with {\tt Macaulay2} \cite{M2} in Examples \ref{ex:18}, \ref{ex:19} and \ref{ex:20}.
We invite our readers to peek at Example \ref{ex:17} where the resultant  from Example \ref{ex:badintroex2} is revisited.

\section{How To Be Bad}
\label{sec2}

We are interested in the subset of the real Grassmannian ${\rm Gr}(k,\mathbb{S}^n)$ whose points are the bad linear spaces $\mathcal{L}$.
By definition, a space $\mathcal{L}$ is bad if  the image cone $\pi_\mathcal{L}(\mathbb{S}^n_+)$ is not closed in $\RR^k$.

\begin{proposition} \label{prop:motivation}
The set of bad linear spaces $\mathcal{L}$ is a semialgebraic 
subset of the Grassmannian ${\rm Gr}(k,\mathbb{S}^n)$.
This set is not closed when $n \geq 3$.
\end{proposition}

\begin{proof}
The first assertion follows from the {\em Tarski-Seidenberg Theorem on Quantifier Elimination} \cite[Theorem A.49]{BPT}.
Indeed, the image $\pi_\mathcal{L}(\SSS^n_+)$ is a semialgebraic subset of $\RR^k$, i.e., it can be described by a Boolean combination of polynomial inequalities. 
The dependence on $\mathcal{L}$ is semialgebraic, as is the statement that the image is not closed. We can eliminate the coordinates of $\RR^k$ to obtain a quantifier-free formula in the Pl\"ucker coordinates of $\mathcal{L}$.
This formula describes the desired semialgebraic subset of the real Grassmannian ${\rm Gr}(k,\mathbb{S}^n)$.

To see that this subset is not closed, fix $k=2,n=3$, let $t$ be a parameter, and consider the quadrics $q_1 = x_1^2$ and $q_2 = x_2^2 + t x_1 x_3$.
Their span is a $2$-dimensional subspace $\mathcal{L}_t$ in $\SSS^3$ for all $t \in \RR$.
For $t \not=0$, the space $\mathcal{L}_t$ is bad because $\pi_{\mathcal{L}_t}(\SSS^3) = \,\{z_1 > 0\} \cup \{z_1 = 0 \,\,{\rm and}\,\, z_2 \geq 0 \}\,$ is not closed.
For $t=0$, the image $\pi_{\mathcal{L}_0}(\SSS^3) =  \RR^2_{\geq 0}$ is closed, so $\mathcal{L}_0$ is good.
This specifies a sequence of bad points in ${\rm Gr}(2,\SSS^3)$ whose limit is a good point, so our set is not closed. 
We note that, by Example \ref{ex:badintroex2}, for $n=2$ the set $\{R = 0\}$ of bad $\mathcal{L}$ is closed in ${\rm Gr}(2,\SSS^3) \simeq \PP^2$.
\end{proof}

A characterization of bad subspaces was given by Pataki. Our first goal is to present~his result.
We consider the spectrahedral cone $\,\mathcal{L} \cap \mathbb{S}^n_+$. 
This convex cone is related to our object of interest 
by duality:

\begin{lemma}\label{dual} The closure of
$ \,\pi_\mathcal{L}(\mathbb{S}^n_+)$ is linearly isomorphic to the cone dual to $\,\mathcal{L} \cap \mathbb{S}^n_+$.
\end{lemma}

\begin{proof}
This follows from the first statement in  \cite[Corollary 16.3.2]{Rock}, 
when $A$ is the linear map $ \pi_{\mathcal L}$ and the convex set $C$ is the PSD cone $\SSS^n_+$.
Note that $\SSS^n_+$ is self-dual.
\end{proof}

We now assume that $k \geq 2$, because $\mathcal{L}$ is always good when $k=1$: the one-dimensional cone $ \pi_\mathcal{L}(\mathbb{S}^n_+)$ equals $ \RR$ if $\,\mathcal{L} \cap \mathbb{S}^n_+ = \{0\}$ and it is $\RR_{\geq 0}$ or $\RR_{\leq 0}$ otherwise.

\begin{definition}
For a linear subspace $\mathcal{L}$ of $\SSS^n$, we fix a quadric $q$ of maximal rank in
the spectrahedral cone $\mathcal{L} \cap \mathbb{S}^n_+$.
The rank of $q$ is an invariant of $\mathcal{L}$, denoted $s= s(\mathcal{L})$ and called the {\em spectrahedral rank} of $\mathcal{L}$.
\end{definition}
If $s = n$ then $\,\mathcal{L} \cap \mathbb{S}^n_+$ is full-dimensional and $\pi_\mathcal{L}(\mathbb{S}^n_+)$ is pointed and closed.
If $s= 0$ then $\,\mathcal{L} \cap \mathbb{S}^n_+ = \{0\}$ and $ \pi_\mathcal{L}(\mathbb{S}^n_+) = \mathcal{L}^\vee $, which is also closed.
Thus we are mostly interested in the cases where $0 < s < n $.
After a linear change of coordinates given by the action of ${\rm GL}(n)$,
we may assume $q = x_1^2 + x_2^2 + \cdots + x_s^2$. 
The matrix that represents a quadratic form $v \in \SSS^n$ has the block structure
$$ \qquad \qquad  V \,=\, \begin{pmatrix} V_{11} & V_{12} \\ V_{12}^T & V_{22} \end{pmatrix}, 
\qquad {\rm where}  \quad V_{11} \in \SSS^s  \,\,\, {\rm and} \,\,\, V_{22} \in \SSS^{n-s}. $$
The following result appears in \cite[Theorems 1,2]{Pataki19}.
A subspace $\mathcal{L}$ is called {\em good} if it is not bad.
We write $\mathcal{L}^\perp = \ker (\pi_{\mathcal{L}})$ for the orthogonal complement of $\mathcal{L}$
in~$\SSS^n$.

\begin{theorem}[Pataki] \label{conditions}
A linear space $\mathcal{L} \subseteq \SSS^n$ is bad  if and only if
there exists a quadric $v \in \mathcal{L}$ whose associated matrix $V$
satisfies $V_{22} \in \SSS^{n-s}_+$ and $\,\im(V_{12}^T) \not\subseteq \im(V_{22})$. 
The space $\mathcal{L}$ is good if and only if there exists 
a positive definite matrix $\,U \in \SSS^{n-s}_+$
 such that $\begin{small} \begin{pmatrix} 0 & 0 \\ 0 & U \end{pmatrix} \end{small} \in \mathcal{L}^\perp$ and,
  for all matrices $V \in \mathcal{L}$, the condition $V_{22} = 0$ implies $V_{12} = 0$.
\end{theorem}

We now present an alternative version of this result.
Since q is positive semidefinite, it is a sum of squares $l_1^2 + \hdots + l_s^2$ of linear forms $l_1, \dots, l_s$, and $I_L$ is the ideal generated by $l_1, \hdots, l_s$.
For instance, if $n \geq 3$ and $q = (x_1-x_2)^2 + (x_1 + x_2)^2 $
then $I_\mathcal{L}  = \langle x_1,x_2 \rangle $.
We consider the inclusions of linear spaces
$$ \mathcal{L} \cap I_\mathcal{L}^2    \,  \,\, \subseteq \,\,\, \mathcal{L} \cap I_\mathcal{L} \,\,\, \subseteq \,\,\,
\mathcal{L} \,\,\subseteq \,\, \mathbb{S}^n \,=\, \RR[x_1,\ldots,x_n]_2 . $$
The first space $\mathcal{L} \cap I_\mathcal{L}^2$ is the linear span of the spectrahedon
$\,\mathcal{L} \cap \mathbb{S}^n_+$, while the second space  $\mathcal{L} \cap I_\mathcal{L}$ 
also records tangent directions relative to the PSD cone. To illustrate the inclusions, we 
consider the bad plane in Example~\ref{ex:badintroex1}, where
$q = x_1^2$, $I_\mathcal{L} = \langle x_1 \rangle $
and $v = x_1 x_2 \in I_\mathcal{L} \backslash I_\mathcal{L}^2$.
We already know that the existence of such a pair $(q,v)$ characterizes non-closed projections.

\begin{corollary} \label{cor:bad} A linear space of quadrics $\mathcal{L} \subset \SSS^n$  is bad if and only if
$s(\mathcal{L}) \geq 1 $ and
\begin{equation}
\label{eq:2cond}  \quad s(\mathcal{L}) + s(\mathcal{L}^\perp)\, < \, n \quad {\rm or} \quad
  \mathcal{L} \cap I_\mathcal{L}^2   \, \subsetneq \,\mathcal{L} \cap I_\mathcal{L}. 
\end{equation}  
\end{corollary}

\begin{proof} We set $s = s(\mathcal{L})$ and assume $I_{\mathcal{L}} = \langle x_1, \dots, x_s \rangle$.
If $s=0$ then $\,\mathcal{L} \cap \mathbb{S}^n_+\,= \,\{0\}$ and thus, by Lemma \ref{dual},
the closure of  $ \pi_\mathcal{L}(\mathbb{S}^n_+)$ is equal to $ \{0\}^\perp = \mathcal{L}^\vee \simeq \RR^k$.
 But, since $\pi_\mathcal{L}(\mathbb{S}^n_+)$ is convex, this implies that
$\pi_\mathcal{L}(\mathbb{S}^n_+) $ equals $ \mathcal{L}^\vee$, so $\mathcal{L}$ is good. 
  Thus, we can now assume $s \geq 1$.

We claim that the two conditions in the disjunction in (\ref{eq:2cond}) are the negations of the two conditions
in the conjunction that characterizes goodness in the last statement of Theorem~\ref{conditions}.
Indeed, since two matrices in $\SSS^n_+$ have trace inner product equal to zero if and only if
their matrix product is the zero matrix, we always have $s(\mathcal{L}^\perp) \leq n-s$. 
The equality $s(\mathcal{L}^\perp) = n-s$ holds if and only if
there is a positive definite $(n-s) \times (n-s)$ matrix $U$  as in Theorem~\ref{conditions}.
Next, consider any matrix $V$ that would correspond to a quadratic form in
$ (\mathcal{L} \cap I_{\mathcal{L}}) \backslash ( \mathcal{L} \cap I_{\mathcal{L}}^2)$.
Containment in $ \mathcal{L} \cap I_{\mathcal{L}}$ means that $V_{22} = 0$,
and non-containment in $\mathcal{L} \cap I_{\mathcal{L}}^2$ means that $V_{12} \not= 0$.
We conclude that Corollary~\ref{cor:bad} is the contrapositive of the last statement of
Theorem~\ref{conditions}.
\end{proof}

\begin{example} \label{ex:threequadrics}
Let $n=k=3$ and let $\mathcal{L}$ be spanned by the span of the quadratic forms
$$ \begin{matrix}
q_1 & = & -52 x_1^2\,+\,412 x_1 x_2+472 x_1 x_3+462 x_2^2+1164 x_2 x_3+750 x_3^2\,, \\
q_2 & = & -101 x_1^2+435 x_1 x_2+480 x_1 x_3+518 x_2^2+1307 x_2 x_3+853 x_3^2\, , \\
q_3 & = &  -55 x_1^2\,+\,362 x_1 x_2+482 x_1 x_3+434 x_2^2+ 1166 x_2 x_3+772 x_3^2\, .
\end{matrix}
$$
These dense quadrics are chosen to hide
the properties of $\mathcal{L}$. They
are identified with symmetric $3 \times 3$-matrices $A_1,A_2,A_3$,
so that $q_i = (x_1,x_2,x_3) A_i (x_1,x_2,x_3)^T$ for $i=1,2,3$.
To reveal the nature of $\mathcal{L}$, we display the elements
$$ \begin{matrix} 7  q_1 -  4  q_2 - 2 q_3  & =  &  6 \cdot (5 x_1 + 7 x_2 + 7 x_3)^2, \qquad \qquad \qquad \qquad \\
17  q_1 -  14  q_2 + 2 q_3  \,\,& =  & \quad 42 \cdot (5 x_1 + 7 x_2 + 7 x_3) (2 x_1 + 5 x_2 + 8 x_3) .\,\,
\end{matrix} $$
Hence $I_\mathcal{L} = \langle 5 x_1 + 7 x_2 + 7 x_3 \rangle$. The space $\,\mathcal{L}  \cap I_\mathcal{L}^2 \,$ is one-dimensional and spanned by the square above, but $ \,\mathcal{L}  \cap I_\mathcal{L} \,$ is two-dimensional. Our two linear combinations form a basis.
The theorem shows that the cone $\pi_\mathcal{L}(\mathbb{S}^3_+)$ is not closed.
It is the union of an open half-space in $3$-space, together with a line through the  origin in the plane that bounds the half-space.
\end{example}

\begin{example}
Let $n=k=3$ and let $\mathcal{L}$ be the span of the quadratic forms
\[q_1 = x_1^2,
 q_2 = x_2^2 + 2x_1x_3,
 q_3 = 2x_2x_3
\]
Hence $I_\mathcal{L} = \langle x_1 \rangle$. 
The space $\mathcal{L} \cap I_\mathcal{L}^2 = \mathcal{L} \cap I_\mathcal{L}$ is spanned by $q_1$. 
The orthogonal complement $\mathcal{L}^\perp$ is spanned by
\[q_1' = 2x_1x_2,
    q_2' = 2x_1x_3 - 2x_2^2,
    q_3' = x_3^2
\]
Hence $I_\mathcal{L} = \langle x_3 \rangle$, so $s(\mathcal{L}) = s(\mathcal{L}^\perp) = 1$ and their sum is $2 < n=3$.
The cone $\pi_\mathcal{L}(\SSS^3_+)$ is the union of an open half-space and a half-line that bounds it.
\end{example}

We now turn the section title around and we focus on {\em how to be good}.
For our best case scenario, we assume $s(\mathcal{L}) =n$.
This means that $\mathcal{L}$ intersects the interior of the
PSD cone $\SSS^n_+$. Hence the intersection $ \mathcal{L} \cap \SSS^n_+$
is a full-dimensional pointed cone in $\mathcal{L} \simeq \RR^k$.
Its convex dual $(\mathcal{L} \cap \SSS^n)^\vee$ is
a full-dimensional pointed cone in 
the dual space $\mathcal{L}^\vee \simeq \RR^k$. 
Lemma \ref{dual} and Theorem \ref{conditions} imply that this closed dual cone is precisely
our projection of the PSD cone:
\begin{equation}
\label{eq:niceimage}
\pi_\mathcal{L} (\SSS^n_+) \,\,\, = \,\,\,
(\mathcal{L} \cap \SSS^n_+)^\vee .
\end{equation} 

Suppose that  $\mathcal{L}$ is generic among points
in ${\rm Gr}(k,\SSS^n)$ that satisfy $s(\mathcal{L}) = n$.
The image cone (\ref{eq:niceimage})  is a generic
spectrahedral shadow, in the sense of \cite{SS}. The
boundary of the image is an algebraic hypersurface that can have 
multiple irreducible components, one for each
matrix rank $r$ in the {\em Pataki range}; see e.g.~\cite[Lemma 5]{FS} and
\cite[Theorem 1.1]{SS}. The degree of the
rank $r$ component is a positive integer, denoted
$\delta(k,n,r)$, that is known as the 
{\em algebraic degree of semidefinite programming}.
These degrees play a major role 
in our main result, which is Theorem \ref{thm:BadCoiso}.
We refer to \cite[Table 1]{SS} for explicit numbers,
and to the bibliographies of \cite{CHS}, \cite{FS}, \cite{NRS}, \cite{SS} and \cite{SU} for additional references.
For instance, for projections of the PSD cone $\SSS^3_+$ into dimensions
$k=3$ and $k=4$, we find that
 $\delta(2,3,2) = 6$ and
$\delta(3,3,r) = 4$ for $r=1,2$.
The following two subspaces exhibit the generic good behavior.

\begin{figure}[h]
$$
        \includegraphics[scale=0.08]{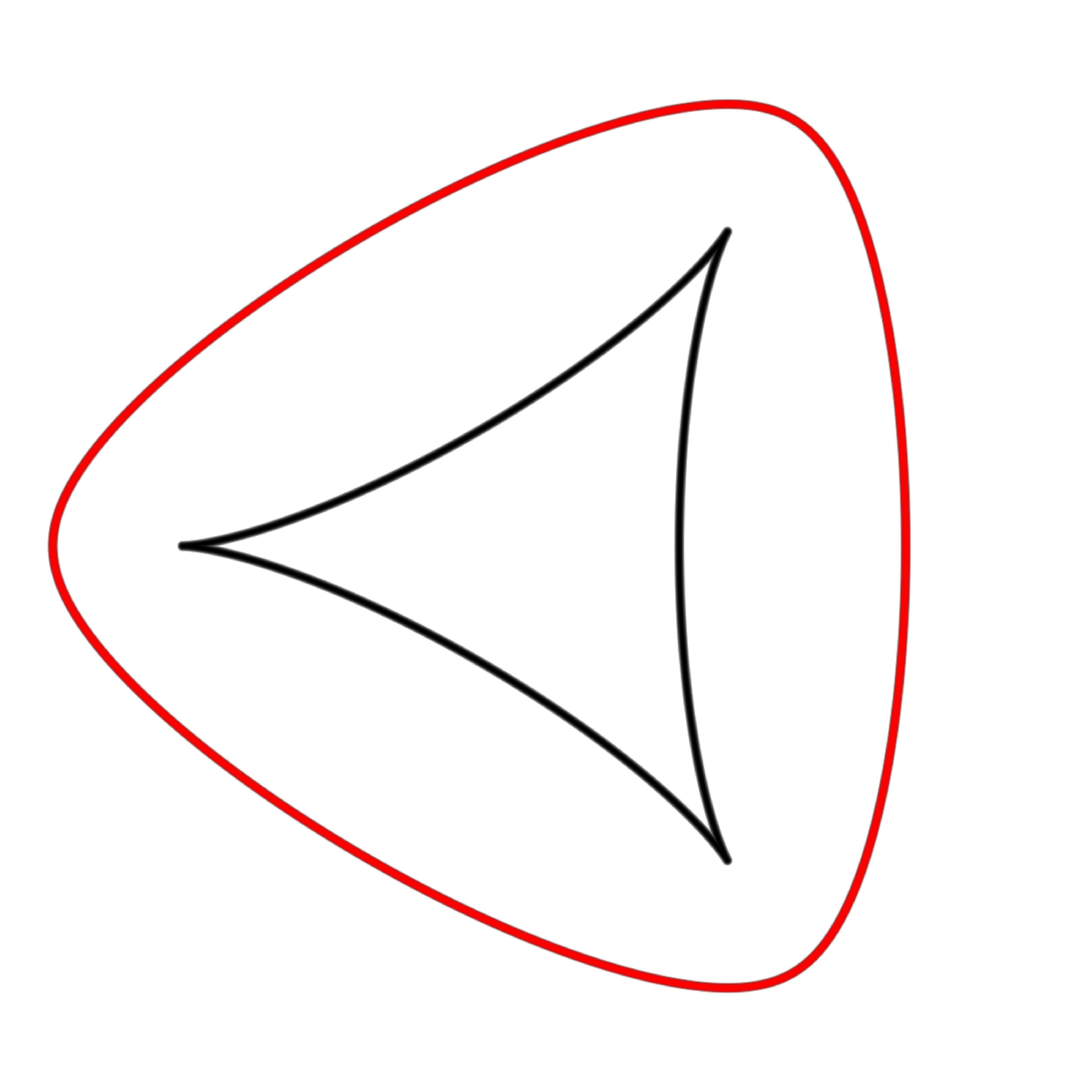} \qquad \qquad \quad
        \includegraphics[scale=0.35]{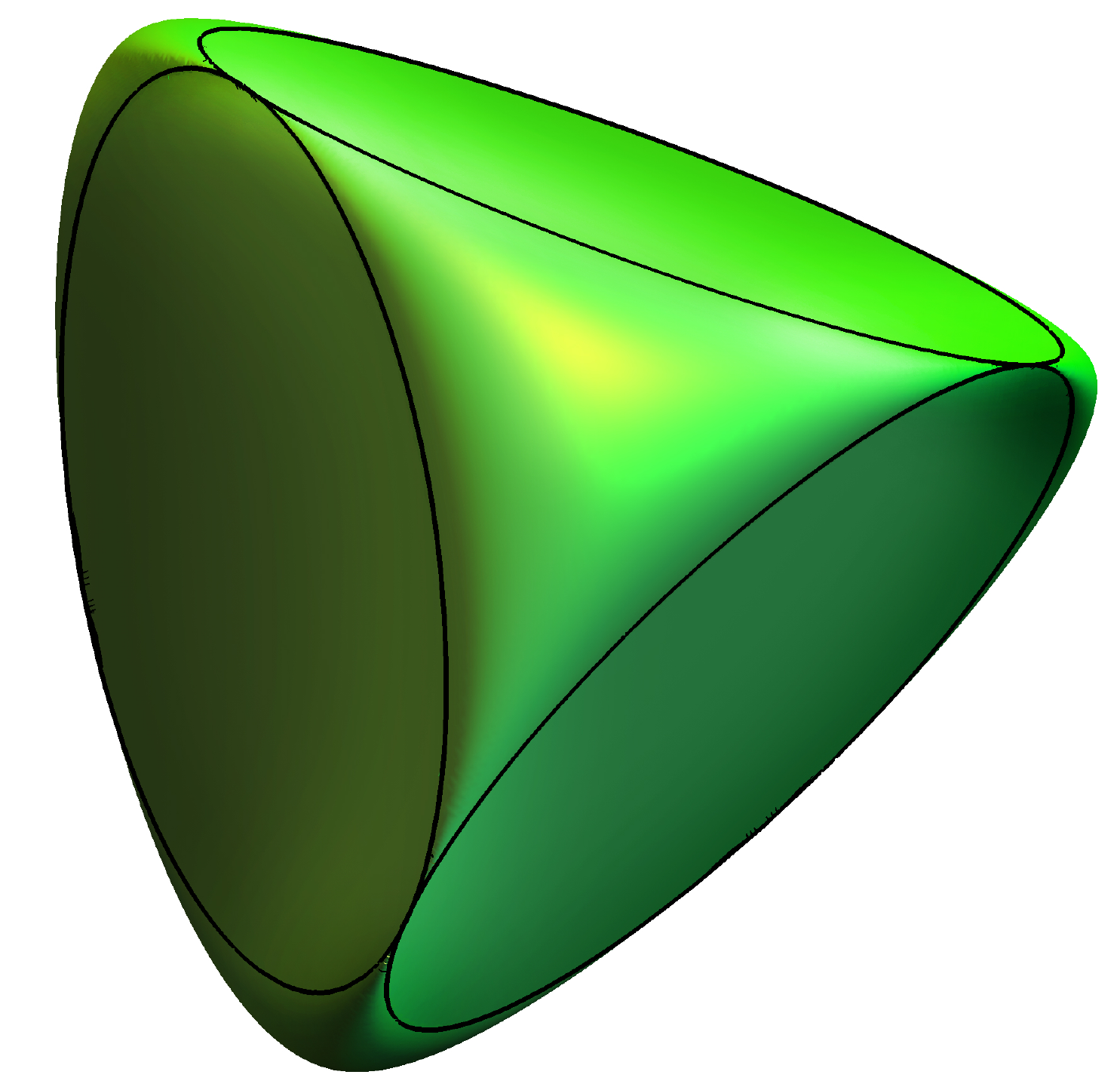}
$$
        \caption{
        The cones over the convex bodies shown
        on the left $(k=3)$ and right $(k=4)$ are
        the images of the $6$-dimensional cone
         $\SSS^3_+$  under  projections $\pi_\mathcal{L}$ defined by
         good subspaces $\mathcal{L}$.
}
        \label{fig:oldimages}
\end{figure}

\begin{example}[$n=k=3$] \label{ex:cone33} Let $\mathcal{L}\,$ be the space
spanned by the rank two quadrics
$$ 
q_1 \,=\, x_1^2 \,+\, (x_2+x_3)^2 \, ,\quad
q_2 \,=\, x_2^2 \,+\, (x_1+x_3)^2 \, , \quad
q_3 \,=\, x_3^2 \,+\, (x_1+x_2)^2 .
$$
Their sum is positive definite, so $s(\mathcal{L}) = 3$.
This specific linear space $\mathcal{L}$ appeared in
\cite{SU} as an illustration for
{\em linear concentration models} in statistics.
In that application,  the three-dimensional cone (\ref{eq:niceimage}) 
serves as the  {\em cone of sufficient statistics} of the model.
Its boundary is an irreducible surface of degree six, defined by 
the polynomial $H_\mathcal{L}$ shown in \cite[Example~1.1]{SU}.
That surface is the cone over the plane sextic curve shown in red
on the left in Figure \ref{fig:oldimages}. For a discussion of this
curve in the context of semidefinite programming see \cite[Example 12.5]{MSbook}.
\end{example}

\begin{example}[$n=3,k=4$] \label{ex:cone34} Let $\mathcal{L}\,$ be the space
spanned by the  rank one quadrics
$$
 q_1 = x_1^2 \, ,\quad
 q_2 = x_2^2 \, ,\quad
 q_3= x_3^2 \, ,\quad
 q_4 = (x_1+x_2+x_3)^2 .
$$
Their sum is positive definite, so $s(\mathcal{L}) = 3$.
The $4$-dimensional cone 
$\mathcal{L} \cap \SSS^3_+$
is the cone over a $3$-dimensional convex body 
known as {\em elliptope} and shown in \cite[Figure 1.1]{MSbook}.
The cone (\ref{eq:niceimage})  is bounded by
four hyperplanes and a quartic threefold in $\RR^4$. These
correspond respectively to the 
four circles and the Roman surface
that bounds the green body on the right in Figure~\ref{fig:oldimages}. 
\end{example}

In the literature, there has been a discrepancy between studies in convex algebraic geometry, like \cite{CHS}, \cite{SS}, and how semidefinite programming is actually used.
The former has focused on generic figures, while the latter often concerns special instances from combinatorial optimization. 
For such scenarios, strong duality can fail, thus motivating works like \cite{KS}, \cite{LP15} and \cite{Pataki19}.
The present paper aims to reconcile these perspectives.
We study special scenarios through the lens
of algebraic geometry, by highlighting 
cases that are generic among the bad ones.

We close this section with a brief exploration of another
connection to algebraic geometry. 
Namely, we consider the restriction of $\pi_{\mathcal{L}}$ to the set of rank one matrices in $\SSS^n_+$. 
This set comprises the extreme rays of the cone $\SSS^n_+$, and it coincides with the image of the map $[\mathcal{L}]_\RR:\RR^n \rightarrow \SSS^n_+$
that takes a vector ${\bf y}$ to the rank one matrix ${\bf y}^T \cdot {\bf y}$.
It is therefore equivalent to study the quadratic map $\RR^n \rightarrow \RR^k$ defined by evaluating the quadrics $q_1,\ldots,q_k$ that span $\mathcal{L}$. 
The image ${\rm im}([\mathcal{L}]_\RR)$ of this map in $\RR^k$ coincides with the image of rank one matrices under~$\pi_{\mathcal{L}}$.
What is this image, and under what conditions on $\mathcal{L}$ is it closed? 
To answer these questions for a small instance, one can apply the method of {\em Cylindrical Algebraic Decomposition} \cite{Col}.
In our experiments we used the implementation {\tt Resolve} in {\tt Mathematica}.

For an algebraic geometer, it is natural to first pass to
the algebraic closure of the field 
$\RR$ and to consider the map
$[\mathcal{L}_\CC]:\CC^n \rightarrow \CC^k$ over the complex numbers.
Our questions remain as above.
What is the image 
 ${\rm im}([\mathcal{L}]_\CC)$ of this map in $\CC^k$,
 and under which conditions on $\mathcal{L}$ is it closed?
 We know from \cite[Theorem 4.23]{MSbook} that
 ${\rm im}([\mathcal{L}]_\CC)$  is closed if
 $\mathcal{L}$ has no zeros in $\PP^{n-1}$.
 To answer these questions for any given instance, 
we can apply the algorithm due to
 Harris, Micha{\l}ek and Sert\"oz \cite{HMS}.
    Our experiments used their   implementation in  {\tt Macaulay2}~\cite{M2}.

We now compare the closure property
for the three sets  $\pi_\mathcal{L}(\SSS^n_+)$,
  ${\rm im}([\mathcal{L}]_\RR)$ and    ${\rm im}([\mathcal{L}]_\CC)$.
In each example we display a basis for $\mathcal{L}$.
One easily finds cases when 
all three sets are~closed, like $(x_1^2,x_2^2)$,
or where none is closed, like $(x_1^2,x_1x_2)$.
The following maps are more interesting.

\begin{example}[Disagreements] \label{ex:disagreements}
Fix $n=3$. Four instances are listed by dimension~$k$:
\begin{itemize}
\item[(2)]
 $\mathcal{L} = (x_1^2+x_2^2,x_1x_3)$. Here 
 $\pi_\mathcal{L}(\SSS^3_+) =   {\rm im}([\mathcal{L}]_\RR) = 
 \{z_1 > 0\} \cup \{z_1=z_2=0\}$ is not closed, but
 the complexification $[\mathcal{L}_\CC]:\CC^3 \rightarrow \CC^2$ is onto. In particular,
$  {\rm im}([\mathcal{L}]_\CC) = \CC^2$ is closed.
\item[(3)] $\mathcal{L} = (x_1x_2,x_1x_3,x_2x_3)$. Here
$\pi_\mathcal{L}(\SSS^3_+) = \RR^3$ is closed because the off-diagonal entries
of a positive semidefinite matrix can be arbitrary. However,
$  {\rm im}([\mathcal{L}]_\RR) $ and  $  {\rm im}([\mathcal{L}]_\CC) $ are not closed.
Their images lack all points that have precisely one coordinate zero.
\item[(4)] $\mathcal{L} = (x_1^2+x_2^2+x_3^2,x_1 x_2, x_2 x_3)$.
Here, $\pi_\mathcal{L}(\SSS^3_+) =   {\rm im}([\mathcal{L}]_\RR) $
is a closed pointed convex cone in $\RR^3$. 
The {\tt Macaulay2} code in \cite{HMS} outputs
$\,{\rm im}([\mathcal{L}]_\CC) \,=\,
\CC^3 \setminus \bigl\{(0,\lambda,\lambda\cdot i)\,:\, \lambda \in \CC \backslash \{0\} \bigr\}$.
This is not closed. However,
since $i = \sqrt{-1}$, none of the missing points is real.
\item[(5)] $\mathcal{L} = (x_1^2-x_2^2,x_3^2,x_1x_2,x_1x_3,x_2x_3)$.
Here $\pi_\mathcal{L}(\SSS^3_+)  $ is not closed because
 $\mathcal{L}^\perp = (x_1^2+x_2^2)$ is a semidefinite ray.
 The complex image $  {\rm im}([\mathcal{L}]_\CC) $ is
 the affine threefold in $\CC^5$ defined by
  $ z_2 z_3=z_4 z_5$ and $ z_1 z_2+z_5^2 = z_4^2$,
  and $  {\rm im}([\mathcal{L}]_\RR) $ is its set of real points.
  Both are closed.
 \end{itemize}
\end{example}

\section{Coisotropic Hypersurfaces}
\label{sec3}

In this section we study the projective variety ${\rm Bad}_{k,n}$.
By definition, this is the Zariski closure in the complex Grassmannian ${\rm Gr}(k,\SSS^n)$ of the set of bad subspaces~$\mathcal{L}$.
Our goal is to characterize ${\rm Bad}_{k,n}$ in terms of objects that are familiar to algebraic geometers. 
We focus on subvarieties of Grassmannians called {\em higher associated varieties}
that are introduced in Section 3.2.E of the book on
{\em Discriminants, Resultants and Multidimensional Determinants}
by Gel'fand, Kapranov and Zelevinsky \cite{GKZ}. They were further studied in recent work of Kohn \cite{kohn} who calls them coisotropic varieties. We shall adopt that name and the notation in \cite{kohn}.

We write $\PP(\SSS^n)  \simeq \PP^{\binom{n+1}{2}-1}$ for the space of symmetric $n \times n$ matrices.
Every point $\mathcal{L}$ in ${\rm Gr}(k,\SSS^n)$ represents a $(k-1)$-dimensional linear subspace $\PP \mathcal{L}$ of $\PP(\SSS^n)$.
Let $\mathcal{Y}$ be a variety of codimension $c$ in $\PP(\SSS^n)$.
For a regular point $X \in {\rm Reg}(\mathcal{Y})$, we write $T_X \mathcal{Y}\,$ for the tangent space of $\mathcal{Y}$ at $X$. 
The {\em $i$-th coisotropic variety} of $\mathcal{Y}$ is 
\begin{equation}
\label{eq:ChDef}
\mathrm{Ch}_i(\mathcal{Y}) \,\,:= \,\,\overline{\bigl\{ \,
\mathcal{L} \in \mathrm{Gr}(c+i, {\SSS}^n) \,\mid \,\exists \,X \in \mathrm{Reg}(\mathcal{Y}) \cap 
\PP \mathcal{L}\,:\, \dim(\PP \mathcal{L} \cap T_X \mathcal{Y}) \geq i \, \bigr\}}.
\end{equation}
The bar denotes  Zariski closure in the Grassmannian. 
Note that $T_X \mathcal{Y}$ has codimension $c$ while
$\PP \mathcal{L} $ has dimension $c+i-1$.
The expected dimension of their intersection is $i-1$. Hence ${\rm Ch}_i(\mathcal{Y})$ is a proper subvariety of $\mathrm{Gr}(c+i, \SSS^n)$, and we expect this to be a hypersurface.
Coisotropic varieties are generalizations of Chow forms.
We shall see this in our examples below. The abbreviation
``Ch'' can thus stand both for Chow from and for coisotropic hypersurface.

The following theorem accomplishes the goal
stated above.
We write $\mathcal{X}_s$ for the subvariety of $\PP(\SSS^n)$ given by all
 symmetric $n \times n$-matrices of rank $\leq s$.
 We set $c = \binom{n-s+1}{2} = \mathrm{codim}(\mathcal{X}_s)$.
 
\begin{theorem} \label{thm:BadCoiso}
The bad subvariety $\,{\rm Bad}_{k,n}$ has codimension one in $\,{\rm Gr}(k,\SSS^n)$. 
It is the union of the irreducible coisotropic hypersurfaces $\,{\rm Ch}_{k-c}(\mathcal{X}_s)$, where $s$ runs over integers such that
\begin{equation}
\label{eq:patakirange}
\binom{n-s+1}{2}\, \,<\,\, k \,\,\leq \,\, \binom{n+1}{2} -\binom{s+1}{2} .
\end{equation}
The degree of the  irreducible polynomial in Pl\"ucker coordinates that defines ${\rm Ch}_{k-c}(\mathcal{X}_s)$ is the algebraic degree of semidefinite programming, which is denoted by $\,\delta(k,n,s)$.
\end{theorem}

The proof will be presented later in this section.
Our first task is to make this statement understandable by
defining all ingredients, and then prove it near the end of this section.
We recall (e.g.~from \cite[Chapter 5]{MSbook})
 that the Grassmannian ${\rm Gr}(k,\mathbb{S}^n)$ is embedded,
 via the {\em Pl\"ucker embedding}, in
a projective space of dimension $\binom{\binom{n+1}{2}}{k}-1$.
 The coordinates we use for ${\rm Gr}(k,\SSS^n)$
are described in \cite[Section 2.1]{kohn}.
The {\em primal Stiefel coordinates} are
the entries of matrices $A_1,\ldots,A_k$ that span $\mathcal{L}$.
The dual Stiefel coordinates are
 matrix entries for a basis of $\mathcal{L}^\perp$.
  If we vectorize these basis elements
and write them as the rows of a matrix with $\binom{n+1}{2}$ columns,
then the maximal minors of this matrix are the
{\em Pl\"ucker coordinates}  of $\mathcal{L}$.

Fix any irreducible variety $\mathcal{Y}$
in $\PP(\mathbb{S}^n) \simeq \PP^{\binom{n+1}{2}-1}$.
Its projectively dual variety $\mathcal{Y}^\vee$ parametrizes hyperplanes that are tangent to $\mathcal{Y}$ at 
some regular point. Note that $\mathcal{Y}$ and~$\mathcal{Y}^\vee$ live in the same ambient space  
$\PP(\mathbb{S}^n)$, since
 $\mathbb{S}^n$ is identified with its dual via the trace~inner product.
We set $c = {\rm codim}(\mathcal{Y})$ and $d = {\rm dim}(\mathcal{Y}^\vee)$.
Following \cite{kohn} and (\ref{eq:ChDef}) above, we write ${\rm Ch}_i(\mathcal{Y})$ for
the $i$-th coisotropic variety  of $\mathcal{Y}$.
This is an irreducible subvariety of the Grassmannian ${\rm Gr}(k,\mathbb{S}^n)$,
where $k=c+i$. The points of ${\rm Ch}_i(\mathcal{Y})$
are linear subspaces~$\mathcal{L}$ that have non-transversal
intersection with the tangent space at some point of $\mathcal{Y}$.
Kohn \cite{kohn}~follows the seminal work of Gel'fand, Kapranov and Zelevinsky \cite{GKZ}
in developing a theory of coisotropic varieties.
She proves in \cite[Corollary 6]{kohn} that ${\rm Ch}_i(\mathcal{Y})$
is a hypersurface if and only if $c \leq k \leq d  + 1$.
In that case, the degree of its equation in Pl\"ucker coordinates equals
\begin{equation}
\label{eq:polardegree} \qquad
 {\rm degree}( {\rm Ch}_i(\mathcal{Y})) \,\, = \,\, \delta_i(\mathcal{Y}) \,:=\,
\hbox{the $i$-th polar degree of $\mathcal{Y}$}.
\end{equation}
This appears in \cite[Theorem 9]{kohn}.
The duality formula in \cite[Theorem 20]{kohn} states 
\begin{equation}
\label{eq:kohnduality}
 {\rm Ch}_i(\mathcal{Y}) \,\simeq \, {\rm Ch}_{d-c+1-i}(\mathcal{Y}^\vee). 
 \end{equation}
This isomorphism is equality if we identify
${\rm Gr}(k,\SSS^n)$ with ${\rm Gr}\bigl(\binom{n+1}{2}-k,\SSS^n\bigr)$
given by  $\mathcal{L} \mapsto \mathcal{L}^\perp$. 
In addition, the polar degree $\delta_i(\mathcal{X})$ is nonzero if and only if $i \leq d - c + 1$.

We now apply these considerations to the
determinantal variety $\mathcal{Y}=\mathcal{X}_s$, the codimension of which equals $c = \binom{n-s+1}{2}$.
It is known that $(\mathcal{X}_s)^\vee = \mathcal{X}_{n-s}$ by
\cite[Proposition~I.1.4.11]{GKZ} or
\cite[Proposition~12]{NRS}.
Hence $d = {\rm dim}(\mathcal{X}_s^\vee) = {\rm dim}(\mathcal{X}_{n-s})
= \binom{n+1}{2}-1 - \binom{s+1}{2}$.
Hence,  by \cite[Corollary 6]{kohn}, ${\rm Ch}_k(\mathcal{X}_s)$ is a hypersurface if and only if
$\binom{n-s+1}{2}\, \,\leq \,\, k \,\,\leq \,\, \binom{n+1}{2} -\binom{s+1}{2} $.
This is almost the same as (\ref{eq:patakirange}), which is known 
as the {\em Pataki range} \cite{FS,NRS,SS}.
However, the minimal value for $k$  in the Pataki range is disallowed in Theorem~\ref{thm:BadCoiso}.
The coisotropic hypersurface for that minimal value of $k$ is the Chow form ${\rm Ch}_0 (\mathcal{X}_s) $.
This does not  contribute to ${\rm Bad}_{k, n}$.

Before proving Theorem \ref{thm:BadCoiso}, let us
explore its implications for matrices of small size.
For $n=2$, the only interesting case is $k=2$.
This was studied in Example 
\ref{ex:badintroex2}, where the resultant $R$ was written
in primal Stiefel coordinates on ${\rm Gr}(2,\SSS^2) = \PP^2$.
The dual Stiefel coordinates are the usual coordinates $(y_0:y_1:y_2)$ on $\PP^2$,
which here agree with the Pl\"ucker coordinates:
$$ y_0 = 2 (a_{12} b_{22} -  a_{22} b_{12}) \, , \quad
y_1 = a_{22} b_{11}  -a_{11} b_{22} \, ,\quad
y_2 = 2 (a_{11} b_{12} -  a_{12} b_{11}). $$
The Veronese curve in $\PP^2$ with equation
$R = y_0 y_2 - y_1^2 $ equals
$\, {\rm Ch}_0(\mathcal{X}_1^\vee) \simeq {\rm Ch}_1(\mathcal{X}_1)$.

\begin{example}[$n=3$] \label{ex:nn33}
We discuss the bad varieties for $k=2,3,4,5$.
Theorem \ref{thm:BadCoiso} states that ${\rm Bad}_{k,3}$ is irreducible and equal 
to the hypersurface ${\rm Ch}_{k-c}(\mathcal{X}_s)$ in ${\rm Gr}(k,\SSS^3) = {\rm Gr}(k,6)$.
The inequalities (\ref{eq:patakirange}) imply
that $s=2$ and $c=1$ for $k=1,2$,
and $s=1$ and $c=3$ for $k=3,4$.

The hypersurface ${\rm Bad}_{2,3}$ has degree $6$ in the $8$-dimensional Grassmannian
 ${\rm Gr}(2,\SSS^3) \subset \PP^{14}$. Its equation is the classical {\em tact invariant} 
of two ternary quadrics $q_1$ and~$q_2$. The tact invariant vanishes
if and only if the conics $\{q_1=0\}$ and $\{q_2=0\}$ are tangent in $\PP^2$. When written in 
the $12=6+6$ entries of the matrices $A_1$ and $A_2$, the tact invariant is a sum of $3210$ 
terms of total degree $12$. This is
the {\em Hurwitz form} of the Veronese surface in $\PP(\SSS^3) = \PP^5$.
The formula in Pl\"ucker coordinates has degree six,
and it appears explicitly in \cite[Example 2.7]{Stu}.

The hypersurface ${\rm Bad}_{3,3} = {\rm Ch}_2(\mathcal{X}_2) \simeq {\rm Ch}_0(\mathcal{X}_1)$ has degree $4$ in 
 ${\rm Gr}(3,\SSS^3) \subset \PP^{19}$. It is the {\em Chow form} of the
 Veronese embedding of $\PP^2$ into $\PP^5= \PP(\mathbb{S}^3)$. This is
 the resultant of three ternary quadrics $q_1,q_2,q_3$ that span~$\mathcal{L}$.
 When written in terms of their $18=6+6+6$ coefficients, this resultant has
 $ 21894 $ terms. We invite our readers to check that this resultant vanishes at the 
 specific instance   $(q_1,q_2,q_3)$ in Example  \ref{ex:threequadrics}. Indeed,
 those three quadrics   have two common zeros in $\PP^2$.
 
 The sextic hypersurface ${\rm Bad}_{4,3} = {\rm Ch}_1 (\mathcal{X}_1)$
 agrees with ${\rm Bad}_{2,3} =  {\rm Ch}_1(\mathcal{X}_2)$ 
 under the identification of ${\rm Gr}(4,\SSS^3)$ with ${\rm Gr}(2,\SSS^3)$.
 If we replace the $15$ Pl\"ucker coordinates
 $p_{ij}$ in \cite[Example 2.7]{Stu} with the complementary
  maximal minors of a $4 \times 6$-matrix,
       then we get an equation of degree $24$ in the 
   entries of a basis  $A_1,A_2,A_3,A_4$ of~$\mathcal{L}$.
 The hypersurface ${\rm Bad}_{5,3} = {\rm Ch}_2(\mathcal{X}_1)$ has degree three in ${\rm Gr}(5,\SSS^3) \simeq  \PP^{5}$. It is 
the determinant hypersurface
$\mathcal{X}_2 = \mathcal{X}_1^\vee$ itself. Indeed, a $5$-dimensional subspace $\mathcal{L}$
is bad if and only if the line $\mathcal{L}^\perp$ is spanned
by a positive semidefinite matrix of rank $\leq 2$.
\end{example}

\begin{example}[$n=4$] \label{ex:vier}
The bad hypersurfaces are irreducible for $n \leq 3$.
The smallest reducible cases arise for $n=4$, with $k=4$ and $k=7$. Namely, we find
$${\rm Bad}_{4,4} \,= \,{\rm Ch}_3(\mathcal{X}_3) \,\cup\, {\rm Ch}_1(\mathcal{X}_2) \quad {\rm and} \quad
{\rm Bad}_{7,4} \,= \, {\rm Ch}_4(\mathcal{X}_2)  \,\cup \,{\rm Ch}_1(\mathcal{X}_1) .$$
We discuss these four coisotropic hypersurfaces and their dual interpretations.
\begin{itemize}
\item ${\rm Ch}_3(\mathcal{X}_3) \simeq {\rm Ch}_0(\mathcal{X}_1)$ has degree $8$ in ${\rm Gr}(4,\SSS^4)$.
The subspace $\mathcal{L} = \RR \{q_1,q_2,q_3,q_4\}$  has a common zero $\mathbf{y}$ in $\PP^3$, i.e.~$\mathcal{L}^\perp$ contains a rank one matrix, namely
${\bf y}^T\cdot {\bf y}$. For a generic bad instance $\mathcal{L}$ in this family, a basis of $\mathcal{L}^\perp$ is given by
\begin{align*}
    q_1 &= x_1^2 + x_2^2 + x_3^2, \\
    q_2 &= \ell_1 x_1 + \ell_2 x_2 + \ell_3 x_3, \\
    q_3 &= \ell_4 x_1 + \ell_5 x_2 + \ell_6 x_3, \\
    q_4 &= \ell_7 x_1 + \ell_8 x_2 + \ell_9 x_3.
\end{align*}
In this example, $\mathbf{y} = (0:0:0:1)$ and $\ell_1,\ell_2,\ldots,\ell_9$ are  generic linear forms.   
\item ${\rm Ch}_1(\mathcal{X}_2)  \simeq {\rm Ch}_3(\mathcal{X}_2)$ has degree $30$ in ${\rm Gr}(4,\SSS^4)$.
In this case, both the zero-dimensional scheme $\PP \mathcal{L} \cap \mathcal{X}_2$ and 
the surface $\PP \mathcal{L}^\perp \cap \mathcal{X}_2$ are singular.
Writing $\ell_1, \dots, \ell_6$ for generic linear forms,
a bad instance in  this family~is
$$
q_1 = x_1^2 + x_2^2 \, ,\,\,\,
q_2 =  \ell_1 x_1 + \ell_2 x_2 \, , \,\,\,
q_3 =  \ell_3 x_1 + \ell_4 x_2 + x_3^2 - x_4^2 \, , \,\,\,
q_4 =  \ell_5 x_1 + \ell_6 x_2 + x_3 x_4 .
$$
The zero-dimensional scheme $\PP \mathcal{L} \cap \mathcal{X}_2$ 
has length $10$ in  $\,\PP(\SSS^4) \simeq \PP^9$. 
The quadric $q_1$ is a point of multiplicity two in that scheme.
\item ${\rm Ch}_4(\mathcal{X}_2) \simeq {\rm Ch}_0(\mathcal{X}_2)$ has degree $10$ in ${\rm Gr}(7,\SSS^4)$.
Here, the threefold $\PP\mathcal{L} \cap \mathcal{X}_2$  is singular and
$\mathcal{L}^\perp $ contains a rank two matrix. That matrix is a singular point
of the plane quartic curve $\PP \mathcal{L}^\perp \cap \mathcal{X}_3$.
For a generic bad instance $\mathcal L$ of this kind, a basis of $\mathcal{L}^\perp$ is given by
$$
q_1 = x_3^2 + x_4^2 \, ,\,\,\,
q_2 = x_1^2 - x_2^2 + \ell_1 x_1 + \ell_2 x_2  \, ,\,\,\,
q_3 = x_1x_2 + \ell_3 x_1 + \ell_4 x_2.
$$
Here $\ell_1, \ell_2, \ell_3, \ell_4$ are binary linear forms in the two unknowns $x_3$ and $ x_4$.
\item ${\rm Ch}_1(\mathcal{X}_1) \simeq {\rm Ch}_2(\mathcal{X}_3)$ has degree $16$ in ${\rm Gr}(7,\SSS^4)$.
Both the zero-dimensional scheme
$\PP \mathcal{L} \cap \mathcal{X}_1$ and the plane quartic  $\PP \mathcal{L}^\perp \cap \mathcal{X}_3$ are singular.
The singular point is a rank three matrix.
For a generic bad instance $\mathcal L$ of this kind, a basis of $\mathcal{L}^\perp$ is given by
$$
q_1 = x_2^2 + x_3^2 + x_4^2 \, ,\,\,\,
q_2 = \ell_1 x_1 + q_2' \, ,\,\,\,
q_3 = \ell_2 x_1 + q_3',
$$
where $\ell_1, \ell_2$ are linear forms in $x_2, x_3, x_4$, and $q_2', q_3'$ are quadrics in $x_2, x_3, x_4$.
\end{itemize}
We note that the two irreducible components of ${\rm Bad}_{7,4}$
are precisely the irreducible factors ${\bf P}$ and ${\bf M}$ of the {\em Vinnikov discriminant}
of a ternary quartic. This was discussed in \cite[Theorem 7.5]{PSV} 
and in \cite[Remark 22]{FKO}.
The three matrices $A,B,C$ that give the determinantal representation of  quartics in \cite{PSV}
span our space $\mathcal{L}^\perp$.
Fl\o ystad, Kileel and Ottaviani \cite{FKO}
present an explicit formula for ${\rm Ch}_4(\mathcal{X}_2)$
as the Pfaffian of a skewsymmetric $20 \times 20$-matrix.
\end{example}

\begin{proof}[Theorem~\ref{thm:BadCoiso}]\label{pf:BadCoiso}
We already saw above that ${\rm Ch}_{k-c}(\mathcal{X}_s)$ is
a hypersurface whenever (\ref{eq:patakirange}) holds. The degree of this coisotropic
hypersurface is the $i$-th polar degree of $\mathcal{X}_s$ by
(\ref{eq:polardegree}), where $i=k-c$.
That degree equals $\delta(k,n,s)$ by \cite[Theorem~13]{NRS}.

In what follows we prove the assertions about the bad variety ${\rm Bad}_{k,n}$
that are stated in the first two sentences of Theorem~\ref{thm:BadCoiso}.
Let $I_{\mathcal{L}}$ and $I_{\mathcal{L}^\perp}$ be defined as in the paragraph prior to Corollary~\ref{cor:bad}. 
The integers $k, n$ are fixed throughout. We begin by assuming that $s$ is in the range (\ref{eq:patakirange}).
We first show that every bad subspace $\mathcal{L} \in {\rm Gr}(k, \SSS^n)$
with $0 < s(\mathcal{L}^\perp) \leq n-s$ and $0 < s(\mathcal{L}) \leq s$ 
is in the coisotropic hypersurface ${\rm Ch}_{k-c}(\mathcal{X}_s)$.
We shall proceed in two~steps.

\smallskip
\noindent
\underbar{Step 1}: Consider a linear space $\mathcal{L} \subset \SSS^n$ with $s(\mathcal{L}) = s$ and $s(\mathcal{L}^\perp) = n-s$.
Then $\mathcal{L}$ is bad if and only if $\mathcal{L} \cap I_\mathcal{L}^2 \subsetneq \mathcal{L} \cap I_\mathcal{L}$, by Corollary~\ref{cor:bad}. 
Points $\mathcal{L}$ such that $\mathcal{L} \cap I_\mathcal{L}^2 = \mathcal{L} \cap I_\mathcal{L}$ lie in a proper subvariety of the Grassmannian $\mathrm{Gr}(k, \SSS^n)$.
This means that the set of bad instances intersects the following set (\ref{eq: set1}) in a Zariski dense subset:
\begin{equation}\label{eq: set1}
\left\{ \,\mathcal{L} \in \mathrm{Gr}(k, \SSS^n)\, \mid \, s(\mathcal{L}) = s \text{ and } s(\mathcal{L}^\perp) = n-s \right\}.
\end{equation}
We may assume, after a change of coordinates, that $\mathcal{L}$ contains the  diagonal matrix
\[ X \,\, = \,\, \mathrm{diag}(\underbrace{1, \dots, 1}_{s}, 0, \dots, 0)\,\, \in \,\,\mathrm{Reg}(\mathcal{X}_s), \] and that $\mathcal{L}^\perp$ contains the special diagonal matrix
\[ Y \,\, = \,\, \mathrm{diag}(0, \dots, 0, \underbrace{1, \dots, 1}_{n-s}). \]
Let $H = T_X(\mathcal{X}_s)$ denote the tangent space of the determinantal variety $\mathcal{X}_s$ at $X$.
This is the linear space of codimension $c = \binom{n-s+1}{2}$ in $\PP(\SSS^n)$ defined 
by $x_{ij} = 0$ for $i >s$ and $j > s$. 
The matrix $Y$ satisfies $Y \in H^\perp$. Hence,
for a generic subspace $\mathcal{L}$ with $X \in \mathcal{L}$ and $Y \in \mathcal{L}^\perp$,
we have $\dim (\PP\mathcal{L} \cap H) \geq k-c$. This number
 exceeds the expected dimension of $\PP\mathcal{L}' \cap H$ among generic $\mathcal{L}'$ that contain $X$.
This shows that any subspace $\mathcal{L}$ with $s(\mathcal{L}) = s$ and $s(\mathcal{L}^\perp) = n-s$
 is a point in the coisotropic variety
$\mathrm{Ch}_{k-c}(\mathcal{X}_s)$. In particular, a generic element of (\ref{eq: set1}) lies in the following full-dimensional semialgebraic set of  real points in ${\rm Ch}_{k-c}(\mathcal{X}_s)$:
\begin{equation}\label{eq: set2}
\left\{\,\mathcal{L} \in \mathrm{Gr}(k, \SSS^n) \mid \exists \,
X \in \mathrm{Reg}(\mathcal{X}_s) \cap \SSS^n_+ \cap \mathcal{L} \,:\, \dim (\PP\mathcal{L} \cap T_X(\mathcal{X}_s)) \geq k- c \,\right\}.
\end{equation}

\smallskip
\noindent
\underbar{Step 2:} If $0 < s(\mathcal{L}) = r < s$ and $0 < s(\mathcal{L}^\perp) = t \leq n-s$. We may assume that
our subspace $\mathcal{L}$ contains the specific diagonal matrix
\[ X \,\, = \,\, \mathrm{diag}(\underbrace{1, \dots, 1}_{r}, 0, \dots, 0) \,\,\in\,\, \mathrm{Sing}(\mathcal{X}_s). \]
Here ${\rm Sing}$ denotes the singular locus. Also,
$\mathcal{L}^\perp$ contains the diagonal matrix
\[ Y \,\, = \,\, \mathrm{diag}(0, \dots, 0, \underbrace{1, \dots, 1}_t). \]
 For any $\epsilon > 0$, let $\mathcal{L}_\epsilon$ be the span
in $\SSS^n$ of  $X^\perp \cap \mathcal{L}$ and
\[
X_\epsilon \,\,=\,\, \mathrm{diag}(\underbrace{1, \dots, 1}_r, \underbrace{\epsilon, \dots, \epsilon}_{s-r}, 0, \dots, 0).
\]
Then $s(\mathcal{L}_\epsilon) = s$ and $\mathcal{L}_\epsilon \to \mathcal{L}$ as $\epsilon \to 0$. 
Let $H = T_{X_\epsilon}(\mathcal{X}_s)$. By construction, we still have $Y \in \mathcal{L}_\epsilon^\perp \cap H^\perp$. By the proof of Step 1, we conclude that $\mathcal{L}_\epsilon \in \mathrm{Ch}_{k-c}(\mathcal{X}_s)$ for all $\epsilon > 0$. Then $\mathcal{L} \in \mathrm{Ch}_{k-c}(\mathcal{X}_s)$ since varieties are closed.
In particular, we see that the set of bad instances with 
$s(\mathcal{L}) + s(\mathcal{L}^\perp) = n$ is dense in the set of all bad instances, so (\ref{eq: set1}) is Zariski dense in (\ref{eq: set2}).

\smallskip

Next we show that values of $s$ outside the range
(\ref{eq:patakirange}) are covered by those in that range.
Suppose that $s$ does not satisfy (\ref{eq:patakirange}).
We claim that the set of bad instances is in the closure of the set of 
bad subspaces for which $s$ is in the range. Fix an $s$ outside the range 
(\ref{eq:patakirange}). Let $\mathcal{L}$ be a bad subspace such that 
either $s(\mathcal{L}) = s$ or $s(\mathcal{L}^\perp) = n-s$.
There are again two cases.

\smallskip
\noindent
\underbar{Case 1:} If $k \leq  c = \mathrm{codim}(\mathcal{X}_s)$, then $k-c \leq 0$. Then ${\rm Ch}_{k-c}(\mathcal{X}_s)$ is either the Chow form or trivially ${\rm Gr}(k, \SSS^n)$, neither of which characterizes bad subspaces as argued below.
By a similar $\epsilon$-argument as that in Step 2, we may assume $s(\mathcal{L}) = s$ and $s(\mathcal{L}^\perp) = n-s$. Then $\mathcal{L} \cap I_\mathcal{L}^2 \subsetneq \mathcal{L} \cap I_\mathcal{L}$ by Corollary~\ref{cor:bad}.
Therefore, if $X \in {\rm Reg}(\mathcal{X}_s) \cap \SSS^n_+ \cap \mathcal{L}$, then $\dim (\PP\mathcal{L} \cap T_X(\mathcal{X}_s)) = \dim \PP(\mathcal{L} \cap I_{\mathcal{L}}) > \dim \PP(\mathcal{L} \cap I_{\mathcal{L}}^2) \geq 0$.
So,  $\dim (\PP\mathcal{L} \cap T_X(\mathcal{X}_s)) \geq 1$.
The set of bad subspaces is contained in the Zariski closure of the following set whose codimension is greater than 1:
\begin{equation*}
\left\{\,\mathcal{L} \in \mathrm{Gr}(k, \SSS^n) \mid \exists \,
X \in \mathrm{Reg}(\mathcal{X}_s) \cap \mathcal{L} \,:\, \dim (\PP\mathcal{L} \cap T_X(\mathcal{X}_s)) \geq 1 \,\right\}.
\end{equation*}
Fix an integer $s < s' < n$ such that $c' = {\rm codim}(\mathcal{X}_{s'}) < k \leq \dim (\mathcal{X}_{n-s'})$.
We claim that the set of bad instances in $\mathrm{Gr}(k, \SSS^n)$ is contained 
in $\mathrm{Ch}_{k-c'}(\mathcal{X}_{s'})$ where $c' = \binom{n-s'+1}{2}$. 
Consider the matrices
\[
X =  \mathrm{diag}(\underbrace{1, \dots, 1}_s, 0, \dots, 0)
\qquad {\rm and} \qquad
X_\epsilon = \mathrm{diag}(\underbrace{1, \dots, 1}_s, \underbrace{\epsilon, \dots, \epsilon}_{s'-s}, 0, \dots, 0).
\]
We may assume $X \in \mathcal{L}$ and we define $\mathcal{L}_\epsilon$ as above. Thus
we have $\mathcal{L}_\epsilon \to \mathcal{L}$ as $\epsilon \to 0$.
Since $s(\mathcal{L}_\epsilon) = s' < n$, we have $n-s' \geq s(\mathcal{L}_\epsilon^\perp) > 0$ by the paragraph folowing Lemma~\ref{dual}.
By the argument in Step 1, we find that $\mathcal{L}_\epsilon \in {\rm Ch}_{k-c'}(\mathcal{X}_{s'})$
for $\epsilon > 0$,
and therefore $\mathcal{L} \in {\rm Ch}_{k-c'}(\mathcal{X}_{s'})$.

\smallskip

\noindent
\underbar{Case 2:}
Let $k > \binom{n+1}{2} - \binom{s+1}{2}$. Then
$\binom{n+1}{2} - k < \binom{s+1}{2} $, so ${\rm Ch}_{\binom{n+1}{2} - k -\binom{s+1}{2}}(\mathcal{X}_{n-s})$ is trivially ${\rm Gr}(\binom{n+1}{2} - k, \SSS^n)$. This does not characterize the orthogonal complement of bad subspaces.
By Corollary~\ref{cor:bad}, bad subspaces are contained in the Zariski closure of the following set whose codimension in the Grassmannian is greater than one:
\[ \left\{\,\mathcal{L} \in \mathrm{Gr}(k, \SSS^n) \mid \exists \,
Y \in \mathrm{Reg}(\mathcal{X}_{n-s}) \cap \mathcal{L}^\perp \,:\, \dim (\PP\mathcal{L}^\perp \cap T_Y(\mathcal{X}_{n-s})) \geq 0 \,\right\}. \]
Choose $0 < s'' < s$ with ${\rm codim}(\mathcal{X}_{n-s''}) \leq \dim \mathcal{L}^\perp < \dim(\mathcal{X}_{s''})$. 
Applying the same $\epsilon$-argument as in Case 1, we find $\mathcal{L}^\perp \in {\rm Ch}_{\binom{n+1}{2}-\binom{s''+1}{2} - k}(\mathcal{X}_{n-s''})$. 
By duality, this means $\mathcal{L} \in {\rm Ch}_{k-c''}(\mathcal{X}_{s''})$ where $c'' = {\rm codim}(\mathcal{X}_{s''})$. 
We have now shown that all bad subspaces $\mathcal{L}$ with $s(\mathcal{L})$ outside the Pataki range lie in one of the coisotropic hypersurfaces ${\rm Ch}_{k-c}(\mathcal{X}_s)$ where $k$ satisfies (\ref{eq:patakirange}).

\medskip

It remains to be seen that the bad subspaces are Zariski dense in
$\mathrm{Ch}_{k-c}(\mathcal{X}_s) $, provided (\ref{eq:patakirange}) holds. By Step 1, the set of bad subspaces is Zariski dense in (\ref{eq: set1}). Since (\ref{eq: set1}) is Zariski dense in (\ref{eq: set2}), it suffices to show that 
(\ref{eq: set2}) is Zariski dense in  $\mathrm{Ch}_{k-c}(\mathcal{X}_s)$. The subset $\mathrm{Reg}(\mathcal{X}_s) \cap \SSS_n^+$ of positive semidefinite rank $s$ matrices is Zariski dense in the variety  $\mathcal{X}_s$. 
The same holds for the incidence variety of pairs $(X, \mathcal{L})$, where $\mathcal{L}$ is tangent to $\mathcal{X}_s$ at $X$. 
We project this incidence variety, and its Zariski dense subset given by $X \in \SSS^n_+$, into the Grassmannian ${\rm Gr}(k, \SSS^n)$. 
The image of the latter is Zariski dense in the image of the former. 
Hence the bad subspaces form a Zariski dense subset of $\mathrm{Ch}_{k-c}(\mathcal{X}_s)$. This completes the proof of Theorem~\ref{thm:BadCoiso}.
\end{proof}
  
\begin{remark}
Our proof gives rise to an explicit parametrization of generic bad subspaces  $\mathcal{L}$ in $\mathrm{Ch}_{k-c}(\mathcal{X}_s)$ which satisfy $s(\mathcal{L}) + s(\mathcal{L}^\perp) = n$.
We shall present a basis for $\mathcal{L}^\perp$, similar to those 
given for $k=7,n=4$ in Example \ref{ex:vier}.
Namely, we start with 
$q_1 =\sum_{i=s+1}^n x_i^2$.
For  $i \geq 2$ we set $q_i = q_i' + \sum_{j=1}^s \ell_{ij} x_j$, where $q_i'$ is a generic quadratic form in $x_1, \dots, x_s$ of trace zero, and the $\ell_{ij}$ are linear forms in $x_{s+1}, \dots, x_n$. The dimension of $\mathcal{L}^\perp$ is supposed to be $\binom{n+1}{2} - k$. 
If $k$ is within the range (\ref{eq:patakirange}) 
then $\dim (\mathcal{L}^\perp) < \binom{n+1}{2} - \binom{n-s+1}{2} = \binom{s+1}{2} + s(n-s)$.
Hence, there is enough freedom to keep all $q_1,q_2,q_3,\ldots$ linearly independent.
The resulting subspaces $\mathcal{L}$ are generically bad,
and they form the Zariski dense subset~(\ref{eq: set1}) of $\mathrm{Ch}_{k-c}(\mathcal{X}_s) $.
\end{remark}

Our main result, Theorem~\ref{thm:BadCoiso}, identifies the
subvarieties in the Grassmannian that are responsible for 
bad behavior in semidefinite programming (SDP). However, these
are complex projective varieties and hence they are one step removed
from the real figures that are of interest in optimization theory.
We close this section by returning to the real and semidefinite setting.
The argument in the last paragraph in the proof of
Theorem~\ref{thm:BadCoiso} gives rise to the following corollary,
aimed at capturing in precise terms what the typical bad subspaces~are.

\begin{corollary} \label{cor:isreal}
Fix integers $n,k,s$ as in(\ref{eq:patakirange}). A generic real subspace
$\mathcal{L}$ in the component ${\rm Ch}_{k-c}(\mathcal{X}_s) $ of
the bad variety $
{\rm Bad}_{k,n}$ is tangent to $\mathcal{X}_s$ at a unique matrix $X \in \SSS^n$,
and this $X$ is real.
A generic space $\mathcal{L}$ such that the matrix
$X$ is positive semidefinite is~bad.
\end{corollary}

\section{Algebraic Computations}
\label{sec4}

In this section we develop computational tools for
the geometric problem studied in this paper.
Suppose we are given matrices $A_1,A_2,\ldots,A_k$
in $\SSS^n$ whose entries are rational numbers. The most basic decision problem
is to determine whether or not
$\mathcal{L} = \RR \{A_1,A_2,\ldots,A_k\}$ is a bad subspace,
and to find a certificate as in Theorem  \ref{conditions}.
The first step in this decision process is the computation of the
spectrahedral rank. Recall that $s(\mathcal{L})$ is the largest rank
of any matrix $X$ in the spectrahedral cone $\mathcal{L} \cap \SSS^n_+$.
For a generic instance $\mathcal{L}$, we have
$s(\mathcal{L}) = 0$ or $s(\mathcal{L}) = n$.
These are the easy cases, where the image
$\pi_\mathcal{L}(\SSS^n_+)$ is either all of $\RR^k$ or a 
closed pointed cone in $\RR^k$. 
We are interested in the decision boundary between these two~regimes.

An upper bound on $s(\mathcal{L})$ is given by the rank of
any matrix $Y $ in $\mathcal{L}^\perp \cap \SSS^n_+$. Indeed,
we have $s(\mathcal{L}) + s(\mathcal{L}^\perp) \leq n $, with
equality for all good  subspaces $\mathcal{L}$, and for generic bad ones.
This leads us to consider the following system of polynomial equations
in $2 \binom{n+1}{2} $ unknowns, which is derived in 
\cite[Corollary 12.12]{MSbook}.

\begin{equation}
\label{eq:criticalvariety}
X \in \mathcal{L} \quad {\rm and} \quad
Y \in \mathcal{L}^\perp \quad {\rm and} \quad X\cdot Y = 0 .
\end{equation}
The pair $(X,Y)$ represents a point in the product space
$\PP \mathcal{L} \times \PP \mathcal{L}^\perp$ inside
   $\PP(\SSS^n) \times \PP(\SSS^n)$.
We call (\ref{eq:criticalvariety}) the {\em critical equations}, and 
its solution set is the  {\em critical variety} of the subspace~$\mathcal{L}$. See \cite[Section 12.2]{MSbook}
for a textbook introduction.
As is customary in algebraic geometry, we work in the  complex projective setting, with
 $\,\PP \mathcal{L} \simeq \PP^{k-1}\,$ and $\,\PP\mathcal{L^\perp} \simeq \PP^{\binom{n+1}{2}-k-1}$.
Thus $\PP \mathcal{L} \times \PP \mathcal{L}^\perp$ is a
variety of dimension $\binom{n+1}{2} -2$.

The  equations (\ref{eq:criticalvariety}) are
reminiscent of the {\em optimality conditions} for SDP,
in the notation used in \cite[(5)]{CHS},
\cite[(12.14)]{MSbook} and
\cite[(3.4)]{NRS}.
However, there is a crucial distinction.
The optimality conditions are based on a flag of subspaces
$\mathcal{V} \subset \mathcal{U}$ where
${\rm dim}(\mathcal{U}) - {\rm dim}(\mathcal{V}) = 2$.
They have many solutions in $\PP(\SSS^n) \times \PP(\SSS^n) $,
counted by  the algebraic degree  of SDP, as shown 
 in \cite{NRS}. In our setting, the flag $\mathcal{V} \subset \mathcal{U}$  is replaced
 by $\mathcal{L} \subseteq \mathcal{L}$. The equations 
  (\ref{eq:criticalvariety}) are expected to have no solutions.
The critical variety of a generic subspace $\mathcal{L}$
is the empty set in $\PP(\SSS^n) \times \PP(\SSS^n)$.
What we care about are the exceptional $\mathcal{L}$ for which (\ref{eq:criticalvariety})
has a solution.

We first review the meaning of the 
equations $X\cdot Y = 0$. These represent
{\em complementary slackness} in SDP.
The {\em normal cycle} of the PSD cone is the semialgebraic~set
\begin{equation}
\label{eq:normalcycle}
{\rm NC}_n \,\, = \,\,
\bigl\{ (X,Y) \in  (\SSS^n_+)^2 \,: \, X \cdot Y = 0 \bigr\}. 
\end{equation}
The normal cycle represents pairs of points in the cone together
with supporting hyperplanes. In the real affine version seen in
(\ref{eq:normalcycle}), this is a semialgebraic set of 
middle dimension $\binom{n+1}{2}$.
If $X$ ranges over matrices of rank $s$ then
$Y$ ranges over complementary matrices of rank $n-s$.
For an algebraic geometer, it is more natural to consider
the complex projective version. This is a reducible
complex algebraic variety, here referred to as
the {\em projective normal cycle}:
\begin{equation}
\label{eq:projectivenormalcycle}
{\rm PNC}_n \,\, = \,\,
\bigl\{ (X,Y) \in  (\PP(\SSS^n))^2 \,: \, X \cdot Y = 0 \bigr\}. 
\end{equation}

\begin{proposition} \label{prop:conormal}
The projective normal cycle ${\rm PNC}_n$ has
$n-1$ irreducible components, each of dimension $\binom{n+1}{2}-2$.
These are the conormal varieties of
the varieties $\mathcal{X}_s$.
\end{proposition}

\begin{proof}
This is the content of \cite[Proposition I.4.11]{GKZ},
revisited for SDP in  \cite[Proposition 12]{NRS} and \cite[Example 5.15]{BPT}.
The dimension statement is
\cite[Proposition 5.10]{BPT}.
\end{proof}

Fix a generic point  $\mathcal{L}$ in the Grassmannian ${\rm Gr}(k,\SSS^n)$.
For dimension reasons, the intersection
$\,(\PP\mathcal{L} \times \PP\mathcal{L}^\perp ) \,\cap \,{\rm PNC}_n\,$
defined by (\ref{eq:criticalvariety}) will be the empty set in  $(\PP(\SSS^n))^2$.
We are interested in subspaces $\mathcal{L}$ for which that intersection is nonempty.

\begin{theorem} \label{thm:LtimesLperp}
The set  $\mathcal{L}$ such (\ref{eq:criticalvariety}) has a solution $(X,Y)$
contains the bad variety:
\begin{equation}
\label{ex:thm17}
\bigl\{ \,\mathcal{L} \in {\rm Gr}(k, \SSS^n)\,:\, 
\,(\PP \mathcal{L} \times \PP\mathcal{L}^\perp ) \,\cap \,{\rm PNC}_n  \,\not= \, \emptyset 
\,\bigr\} \,\, \supseteq \,\,{\rm Bad}_{k,n}.
 \end{equation}
Equality holds unless $k=\binom{n-s+1}{2}$, when
the difference is the Chow form ${\rm Ch}_0(\mathcal{X}_s)$.
The closure of the set of bad subspaces 
is the set of real $\mathcal{L}$ for which 
$\mathcal{L} \times \mathcal{L}^\perp $ intersects the 
normal cycle ${\rm NC}_n$ nontrivially, i.e.~(\ref{eq:criticalvariety}) has
a solution  $(X,Y)$ with $X {\not=} \,0,Y {\not=} \,0$.
\end{theorem}

\begin{proof}
This follows from Theorem \ref{thm:BadCoiso} and  Corollary \ref{cor:isreal},
using the fact about conormal varieties stated in Proposition \ref{prop:conormal}.
It suffices to show that a Zariski dense subset of 
the bad variety ${\rm Bad}_{k, n}$ is contained in the left hand side of (\ref{ex:thm17}).
We choose that subset to be the union of (\ref{eq: set2}) over all $s$ within the range~(\ref{eq:patakirange}). Fix $s$ and consider $\mathcal{L}$ in (\ref{eq: set2}). Pick $X \in {\rm Reg}(\mathcal{X}_s) \cap \mathcal{L}$ such that $\dim (\PP\mathcal{L} \cap T_X(\mathcal{X}_s)) \geq k-c$. 
Let $I_{\mathcal{L}}$ and $I_{\mathcal{L}^\perp}$ be defined as in the paragraph prior to Corollary~\ref{cor:bad}. 
Let $\mathcal{L}'$ be the image of $\mathcal{L}$ under the map $\SSS^n \to \SSS^{n-s}$ which restricts the domain of each matrix to the variety of $I_{\mathcal{L}}$. 
Then $\dim (\mathcal{L}' )\leq c-1 = \dim(\SSS^{n-s}) -1$,
 so that $\pi_{\mathcal{L}'}: \SSS^{n-s} \to \mathcal{L}'^\vee$ has a nontrivial kernel.
  If $Y \in \ker (\pi_{\mathcal{L}'})$, then $Y \in \mathcal{L}^\perp$ and $X \cdot Y = 0$.

The same relationship between 
coisotropic hypersurfaces and the conormal variety extends from 
our specific varieties $\mathcal{X}_s$ to
 arbitrary projective varieties. This is essentially biduality, and we view it
 as a  geometric refinement of \cite[Section 4]{kohn}.
\end{proof}

\begin{remark}
The  ``unless'' statement in the second sentence of Theorem \ref{thm:LtimesLperp}
looks mysterious at first sight. We here offer an explanation for the case $n=k=3$ and $s=1$.
The left hand side of (\ref{ex:thm17}) equals ${\rm Ch}_2(\mathcal{X}_2) \,\cup\, {\rm Ch}_0(\mathcal{X}_1)$.
We now derive the irreducible polynomials for the two  components. 
Following \cite[Example 5.3]{MSbook},
we introduce nine affine coordinates $a,b,\ldots,h,i$  on ${\rm Gr}(3,\SSS^3)$. 
To this end, we fix bases 
$\{q_1,q_2,q_3\}$ for $\mathcal{L}$ and
$\{r_1,r_2,r_3\}$ for $\mathcal{L}^\perp$ as~follows:
$$ \begin{matrix}
q_1 = x_1^2+ax_1 x_2+b x_1 x_3+c x_2 x_3\,, & \quad &  r_1 = a x_1^2+d x_2^2+g x_3^2-2 x_1 x_2,\\
q_2 = x_2^2+dx_1 x_2+e x_1 x_3+ f x_2 x_3\,, & \quad & r_2 = b x_1^2+e x_2^2+h x_3^2-2 x_1 x_3, \\
q_3 = x_3^2+g x_1 x_2+h x_1 x_3+i x_2 x_3\,, &  \quad &  r_3 = c x_1^2+f x_2^2+i x_3^2 - 2 x_2 x_3. \\
\end{matrix}
$$
The coisotropic hypersurface ${\rm Ch}_2(\mathcal{X}_2)$ is defined by a resultant with $218$ terms:
$$ \begin{small} \begin{matrix} {\rm Res}(q_1,q_2,q_3)  & \! = \!\! &  a^2 b d e f h i^2-a^2 b d f^2 h^2 i-a^2 b e^2 f g i^2+a^2 b e f^2 g h i-a^2 c d
 e^2 h i^2+a^2 c d e f h^2 i \\  & & +\,a^2 c e^3 g i^2- \cdots 
  +a f h+b d i+b f g+c d h+3 c e g-2 a d-2 b h-2 f i+1. 
  \end{matrix} \end{small}
  $$
The coisotropic hypersurface ${\rm Ch}_0(\mathcal{X}_1)$ is defined by a
resultant with $549$ terms:
$$ \begin{small} \begin{matrix} {\rm Res}(r_1,r_2,r_3) &  = & 
a^4 e^4 i^4-4 a^4 e^3 f h i^3+6 a^4 e^2 f^2 h^2 i^2-4 a^4 e f^3 h^3 i+a^4 f^4 h^4
-4 a^3 b d e^3 i^4\\ & & \!\! +\,12 a^3 b d e^2 f h i^3 - \cdots
+384 c e f g i
-128 c f^2 g h+256 a b e g+256 c e h i-512 c e g.
\end{matrix}
\end{small}
$$
The bad variety ${\rm Bad}_{3,3}$ equals  ${\rm Ch}_2(\mathcal{X}_2)$. All generic subspaces
in ${\rm Ch}_0(\mathcal{X}_1)$ are good.

The objects in (\ref{eq:criticalvariety}), (\ref{eq:normalcycle}) and
(\ref{eq:projectivenormalcycle}) are symmetric under switching the two factors.
But this is not the case when it comes to $\pi_\mathcal{L}(\SSS^n_+)$ being closed.
For deciding between bad and good, the symmetry between primal
and dual is broken. From an SDP perspective, one can see this 
in Proposition 1 of \cite[Section 3]{Pataki19}.
For a concrete example, consider the dual pair
$\mathcal{L} = \RR\{x_1^2+x_2^2,x_1 x_2+x_1 x_3,x_2 x_3 \}$
and $\mathcal{L}^\perp = \RR \{x_3^2,x_1^2-x_2^2,x_1x_2-x_1x_3 \}$.
Then $\mathcal{L} \in {\rm Ch}_2(\mathcal{X}_2) \backslash {\rm Ch}_0(\mathcal{X}_1)$ is bad, 
and quite typical for this, whereas
$\mathcal{L}^\perp \in {\rm Ch}_1(\mathcal{X}_0) \backslash  {\rm Ch}_2(\mathcal{X}_2)$ is good.
\end{remark}

We can use Theorem \ref{thm:LtimesLperp} to compute equations that define
our coisotropic hypersurfaces. Namely, consider the incidence variety in 
$\PP(\SSS^n) \times \PP(\SSS^n) \times {\rm Gr}(k,\SSS^n)$ that is 
defined by the critical equations (\ref{eq:criticalvariety}).
The hypersurface we are interested in is the image of that incidence variety under the map
$\,\PP(\SSS^n) \times \PP(\SSS^n) \times {\rm Gr}(k,\SSS^n) \rightarrow {\rm Gr}(k,\SSS^n)$.
In particular, if we fix some particular rank $s$, then the image of 
the incidence variety in
$\,\mathcal{X}_s \times \mathcal{X}_{n-s} \times {\rm Gr}(k,\SSS^n)$ is the
 irreducible hypersurface ${\rm Ch}_{k-c}(\mathcal{X}_s)$ in  $ {\rm Gr}(k,\SSS^n)$.
 Algebraically, one obtains the polynomial defining ${\rm Ch}_{k-c}(\mathcal{X}_s)$
 by eliminating $X$ and $Y$ from the following {\em  rank $s$ critical equations of~$\mathcal{L}$}:
  \begin{equation}
\label{eq:criticalvarietyranks}
X \in \mathcal{L} \,\,\, {\rm and} \,\,\,
Y \in \mathcal{L}^\perp \,\,\, {\rm and} \,\,\, X\cdot Y = 0 \,\,\, {\rm and} 
\,\,\, {\rm rank}(X) \leq s \,\,\, {\rm and} \,\,\, {\rm rank}(Y) \leq n-s. 
\end{equation}

\begin{corollary} \label{cor:twenty}
Under the  assumption of Corollary \ref{cor:isreal},
the system (\ref{eq:criticalvarietyranks}) has
a unique solution $(X,Y)$ in $(\PP(\SSS^n))^2$.
Here $X$ is the point of tangency and $Y$
is in the normal~space.
\end{corollary}
 
\begin{example} \label{ex:17}
This {\tt Macaulay2} code represents the system (\ref{eq:criticalvarietyranks})
for $n=k=2, s=1$:
\begin{small}
\begin{verbatim}
R = QQ[ x1,x2, y11,y12,y22, a11,a12,a22, b11,b12,b22 ];
A = matrix {{a11,a12},{a12,a22}};  B = matrix {{b11,b12},{b12,b22}};
X = x1*A + x2*B;                   Y = matrix {{y11,y12},{y12,y22}};
I = ideal(trace(A*Y),trace(B*Y)) + minors(1,X*Y) + minors(2,X) + minors(2,Y)
\end{verbatim}
\end{small}
The following eliminates the pair $(X,Y)$ and retains the subspace
$\mathcal{L} = {\rm span}({\tt A},{\tt B})$.
\begin{small}
\begin{verbatim}
eliminate({ x1,x2, y11,y12,y22 }, saturate( I, ideal(y11,y12,y22) ) )
\end{verbatim}
\end{small}
As predicted, the output is the resultant $R = {\rm Bad}_{2,2} = {\rm Ch}_1(\mathcal{X}_1)$
from Example~\ref{ex:badintroex2}.
Note the importance of the saturation step in 
representing subschemes of $(\PP(\SSS^n))^2$.
\end{example}

\begin{example}[$n=3,k=4,s=1$] \label{ex:18}
We represent $\mathcal{L}^\perp$ by a
basis of $3 \times 3$-matrices ${\tt U}$ and~${\tt V}$.
Their $6+6$ entries are the dual Stiefel coordinates on
$ {\rm Gr}(4,\SSS^3)$. We run
\begin{small}
\begin{verbatim}
R = QQ[ y1,y2, x1,x2,x3,x4,x5,x6, u1,u2,u3,u4,u5,u6, v1,v2,v3,v4,v5,v6 ]
X = matrix {{x1,x2,x3},{x2,x4,x5},{x3,x5,x6}}
U = matrix {{u1,u2,u3},{u2,u4,u5},{u3,u5,u6}}
V = matrix {{v1,v2,v3},{v2,v4,v5},{v3,v5,v6}}
Y = y1*U + y2*V
I = ideal(trace(X*U),trace(X*V)) + minors(1,X*Y) + minors(2,X) + minors(3,Y)
I = I:minors(1,X);  I = I:ideal(y1,y2);
eliminate({ y1,y2, x1,x2,x3,x4,x5,x6 }, I )
\end{verbatim}
\end{small}
The output  is a polynomial in {\tt u1..u6,v1..v6}
with $3210$ terms of degree $12$.
This is the tact invariant we saw in Example \ref{ex:nn33}. It
 defines ${\rm Bad}_{4,3} = {\rm Ch}_1(\mathcal{X}_1) \simeq {\rm Ch}_1(\mathcal{X}_2)$.
\end{example}

We close with two case studies on numerical examples.
In both cases, the  critical variety defined by (\ref{eq:criticalvarietyranks}) 
consists of a single rational point, and it is computed using the
command {\tt criticalIdeal} in the {\tt Macaulay2} package
{\tt SemidefiniteProgramming} \cite{CKP}. By using this package,
we can compare the algebraic
approach described above with 
Pataki's facial reduction \cite{Pataki13} in 
the usual numerical framework of SDP.

\begin{example}[$n=k=4$] \label{ex:19} \ Let $t$ be a parameter
and $\mathcal{L}_t$ the space with~basis
$$ \begin{tiny}
 \! A_1 = \begin{bmatrix}
 180 \! & \! 112 \! & \! 205 \! & \! 131\\ 112 \! & \! 88 \! & \! 131 \! & \! 96 \\
  205 \! & \! 131 \! & \! 228 \! & \! 152 \\ 131 \! & \! 96 \! & \! 152 \! & \! 104
 \end{bmatrix} \! , \,
   A_2 = \begin{bmatrix}
  428 \! & \! 253 \! & \! 473 \! & \! 288 \\ 253 \! & \! 238 \! & \! 262 \! & \! 227 \\
  473 \! & \! 262 \! & \! 516 \! & \! 307 \\ 288 \! & \! 227 \! & \!  307 \! & \! 168 \\
    \end{bmatrix} \! , \,
 A_3 = \begin{bmatrix} 216 \! & \! 123 \! & \! 234 \! & \! 137 \\
 123 \! & \! 128 \! & \! 118 \! & \! 116 \\ 234 \! & \! 118 \! & \! 252 \! & \! 138 \\ 137 \! & \! 116 \! & \! 138 \! & \! 68 
 \end{bmatrix} \!, \,
 A_4 = \begin{bmatrix} 320 \!& \! t \!&\! 380 \!&\! 254 \\
 t \! & \! 140 \! & \! 258 \! & \! 166 \\  380 \! & \! 258 \! & \! 448 \! & \! 342 \\
  254 \! & \! 166 \! & \! 342 \! & \! 208 \end{bmatrix}\! .
  \end{tiny}
$$
Is there a value of $t$ for which $\mathcal{L}_t$ is bad?
Geometrically, $\{\mathcal{L}_t\}_{t \in \RR}$ is a line in ${\rm Gr}(4,\SSS^4) \subset \PP^{209}$.
What is its intersection  with ${\rm Bad}_{4,4}$?
To answer this question, we consider $s=2$ and we  evaluate the Hurwitz form
${\rm Ch}_1(\mathcal{X}_2)$ on $\mathcal{L}_t$; see \cite{Stu}. The result is 
a primitive polynomial $f \in \ZZ[t]$ of degree $30$.
One of its roots is $194$.
Each coefficient has over $100$ digits. The leading coefficient equals
\begin{small}
34006196837917896573795931\end{small}\begin{small}
713719797442228459580467476929073669732608826214924094017413181247522771810155185.
\end{small}

We now substitute $t=194$,   and we continue our computation in
   the polynomial ring $\QQ[x_1,x_2,x_3,x_4,y_{11}, y_{12}, \ldots,y_{44}]$.
  We write  $X = x_1 A_1 + x_2 A_2 + x_3 X_3 + x_4 X_4$ for a general
 matrix in $\mathcal{L}$. The 
 constraint $Y \in \mathcal{L}^\perp$ is encoded by the four linear equations 
 $ A_i \circ Y  = 0$
 in the ten entries of the matrix $Y = (y_{ij})$.
The ideal for (\ref{eq:criticalvarietyranks})
 is given by the $3\times 3$-minors of $X$ and $Y$
 as well as the   entries of $X\cdot Y $.
After saturating by the irrelevant ideal  of $\PP\mathcal{L} \times \PP(\SSS^4)$,
we obtain the homogeneous maximal ideal of the unique point $(X^*,Y^*)$ in the critical variety:
$$\begin{smallmatrix} \langle \, x_1-14x_4, x_2+18x_4, x_3-24x_4\,,  \,\,
 314 y_{11} -197 y_{44} ,314 y_{12}-11 y_{44}, 314 y_{12}-11 y_{44}, 314y_{13}+213 y_{44},
 \\ \,\,\,\,314 y_{14}-53 y_{44},157 y_{22}-313 y_{44},  157 y_{23} +30 y_{44},157 y_{24} +215 y_{44}, 157 y_{33}-117 y_{44}, 157 y_{34} +12 y_{44} \, \rangle.
\end{smallmatrix}$$
The matrices $X^*$ and $Y^*$ are unique up to scaling. They are positive semidefinite of rank~$2$.
For instance, $X^* = 14 A_1 - 18 A_2 + 24 A_3 + A_4$ is the point of tangency for $\mathcal{L}$
at the variety~$\mathcal{X}_2$.
For an alternative view, we consider the degenerate SDP
\begin{equation}
\label{eq:minimize}
{\rm Minimize}\,\, 0 \circ Y \quad \hbox{subject to} \,\,\,\,  A_1 \circ Y =  A_2 \circ Y =A_3 \circ Y =A_4 \circ Y = 0 . 
\end{equation}
We can solve this numerically in {\tt Macaulay2}, using the following commands:
\begin{small}
\begin{verbatim}
needsPackage "SemidefiniteProgramming"
R = QQ[x1, x2, x3, x4]
A = matrix{
{180*x1+428*x2+216*x3+320*x4,112*x1+253*x2+123*x3+194*x4,
        205*x1+473*x2+234*x3+380*x4,131*x1+288*x2+137*x3+254*x4},
{112*x1+253*x2+123*x3+194*x4, 88*x1+238*x2+128*x3+140*x4,
        131*x1+262*x2+118*x3+258*x4, 96*x1+227*x2+116*x3+166*x4},
{205*x1+473*x2+234*x3+380*x4,131*x1+262*x2+118*x3+258*x4,
        228*x1+516*x2+252*x3+448*x4,152*x1+307*x2+138*x3+342*x4},
{131*x1+288*x2+137*x3+254*x4, 96*x1+227*x2+116*x3+166*x4,
        152*x1+307*x2+138*x3+342*x4, 104*x1+168*x2+68*x3+208*x4}}
objFun = 0*x1 + 0*x2 + 0*x3 + 0*x4
P = sdp({x1, x2, x3, x4}, A, objFun)
(Y, x, X, v) = optimize P
\end{verbatim}
\end{small}
The numerical output approximates the point $(X^*,Y^*)$  in the critical variety of $\mathcal{L}$.
\end{example}

\begin{example}[$n= 4, k = 5, s=2 $] \label{ex:20}
Let $A_1,\ldots,A_5$ be the five matrices displayed in 
equation (2.1) of \cite[Example 4]{Pataki19}. Their linear span
$\mathcal{L} \in {\rm Gr}(5,\SSS^4) $ is a bad subspace.
However,  $\mathcal{L}$ is not generic in the sense of
Corollary \ref{cor:isreal}. To see this, we compute
the saturation of the ideal specified in (\ref{eq:criticalvarietyranks}).
We find that the critical variety is a reducible surface in 
$\PP \mathcal{L} \times \PP \mathcal{L}^\perp$.
We conclude that $\mathcal{L}$
is a point in ${\rm Bad}_{5,4}$, but it does not
satisfy the hypothesis in Corollary \ref{cor:twenty}, as that would imply
that the variety is only one point. The projection of the critical variety into
$\PP \mathcal{L} \simeq \PP^4$ is the plane defined by
$\langle  3 x_1+ 2 x_3+3 x_4+9 x_5, 3 x_2-x_3+3 x_5 \rangle $.
The point $x = (-1,-1,0,-2,1)$ lies in that plane. It specifies the 
matrix $X = \sum_{i=1}^5 x_i A_i = {\rm diag}(1,1,0,0)$,
seen on the right in \cite[equation (2.4)]{Pataki19}. Note that $X$
is in $ \mathcal{L} \cap \SSS^4_+$ and 
${\rm rank}(X) = s(\mathcal{L})=2$.

The bad variety ${\rm Bad}_{5,4}$ is the 
hypersurface $  {\rm Ch}_2(\mathcal{X}_2)$, which is
self-dual and has degree $\delta(5,4,2) = 42$.
To construct typical points, we 
choose random  $B_1,B_2,B_3$ in $\SSS^4$,
and we replace $A_1$ by $B_1$ and $A_5$ by $B_2 + t B_3$,
where $t$ is a new unknown. Let $\mathcal{L}_t$ denote the resulting subspace.
We repeat the above computation of (\ref{eq:criticalvarietyranks}) for  $\mathcal{L}_t$.
By eliminating all $15$ variables $x_i$ and $ y_{jk}$, we obtain a principal ideal
$\langle f(t) \rangle $, where $f \in \ZZ[t]$ has degree $42$.
This is the restriction of ${\rm Ch}_2(\mathcal{X}_2)$ to the line 
$\{ \mathcal{L}_t \}_{t \in \CC}$ in $ {\rm Gr}(5,\SSS^4) \subset \PP^{251}$.
The real roots of $f(t)$ are the candidates for bad subspaces $\mathcal{L}_t$.

We experimented with the degenerate SDP 
given by $A_1,A_2,A_3,A_4,A_5$.
Using {\tt SemidefiniteProgramming} in {\tt Macaulay2}  \cite{CKP}, \cite{M2},
we found that the critical ideal has codimension 11 and degree 11. It was faster to compute the projection of the critical variety into $\PP \mathcal{L} \simeq \PP^4$ after the change of coordinates  in \cite[(2.3)]{Pataki19}.
However, the numerical solver performed better before the change of coordinates, and it gave

 \begin{tiny}
\begin{verbatim}
      X = | .638388     -.226356    2.80061e-8  -1.09986e-8 |        Y = | 9.37872e-10 1.08292e-10 -1.80616e-7  .000428954   |
          | -.226356    .0802601    1.70075e-8  -9.00272e-9 |            | 1.08292e-10 6.92506e-10 -1.7803e-7   .000273252   |
          | 2.80061e-8  1.70075e-8  2.06286e-8  -1.99587e-9 |            | -1.80616e-7 -1.7803e-7  .000703099   -.0000774043 |
          | -1.09986e-8 -9.00272e-9 -1.99587e-9 6.2732e-10  |            | .000428954  .000273252  -.0000774043 538.187      |
\end{verbatim}
\end{tiny}
This is the primal-dual pair $(X^*,Y^*)$ of numerical solutions to our critical equations (\ref{eq:criticalvarietyranks}). Rounding small numbers down to zero, we see that approximately ${\rm rank} (X^*) = 2$ while ${\rm rank}(Y^*) = 1$.
This is consistent with theoretical analysis.
If the subspace $\mathcal{L}$ were generic in ${\rm Bad}_{5,4}$,
then the rank of these matrices add up to $n=4$. They would be unique up to scaling and their entries would be algebraic numbers of degree $42 = \delta(5,4,2)$ over the rationals $\QQ$.
\end{example}

\begin{acknowledgements}
We thank Diego Cifuentes, Kathl\'en Kohn, Gabor Pataki
and Raphael Pellegrin for helpful conversations about topics featured in this article.
\end{acknowledgements}

\end{document}